\documentclass[12pt]{amsart}
\usepackage{amscd,amssymb,amsthm,verbatim}
\newcommand{\C}{{\mathbb C}}

\newcommand{\diam}{\operatorname{diam}}
\newcommand{\Diff}{\operatorname{Diff}}

\newcommand{\dvol}{\operatorname{dvol}}

\newcommand{\Id}{\operatorname{Id}}

\newcommand{\inj}{\operatorname{inj}}

\newcommand{\Isom}{\operatorname{Isom}}

\newcommand{\Nil}{\operatorname{Nil}}

\newcommand{\R}{{\mathbb R}}

\newcommand{\Ric}{\operatorname{Ric}}
\newcommand{\Riem}{\operatorname{Riem}}

\newcommand{\SL}{\operatorname{SL}}
\newcommand{\SO}{\operatorname{SO}}
\newcommand{\Sol}{\operatorname{Sol}}

\newcommand{\SU}{\operatorname{SU}}

\newcommand{\Tr}{\operatorname{Tr}}

\newcommand{\Z}{{\mathbb Z}}

\numberwithin{equation}{section}
\theoremstyle{plain}
\newtheorem{definition}[equation]{Definition}

\newtheorem{lemma}[equation]{Lemma}
\newtheorem{theorem}[equation]{Theorem}
\newtheorem{proposition}[equation]{Proposition}
\newtheorem{corollary}[equation]{Corollary}

\theoremstyle{remark}

\newtheorem{remark}[equation]{Remark}
\newtheorem{example}[equation]{Example}
\usepackage{graphicx}
\input{amssym.def}
\input{amssym}
\input{epsf}
\usepackage{a4wide}

\begin{document}

\title[On the long-time behavior of type-III Ricci flow solutions]
{On the long-time behavior of type-III Ricci flow solutions}

\author{John Lott}
\address{Department of Mathematics\\
University of Michigan\\
Ann Arbor, MI  48109-1109\\
USA} \email{lott@umich.edu}

\thanks{This work was
supported by NSF grant DMS-0306242}
\date{March 11, 2007}

\begin{abstract}
We show that  three-dimensional homogeneous
Ricci flow
solutions that admit finite-volume quotients
 have long-time limits given by
expanding solitons.  We show that the same is true for a
large class of four-dimensional homogeneous
solutions. We give an extension of Hamilton's
compactness theorem that does not assume a lower 
injectivity radius bound, in terms of Riemannian groupoids.
Using this,
we show that the long-time behavior of type-III Ricci flow
solutions is governed by the dynamics of an $\R^+$-action on
a compact space.
\end{abstract}

\maketitle

\section{Introduction}

A type-III Ricci flow solution is a $1$-parameter family 
$\{g(t)\}_{t \in (0, \infty)}$ of Riemannian metrics on a manifold $M$
that satisfy the Ricci flow equation and have sectional
curvatures that decay at least as fast as $t^{-1}$, i.e.
$\sup_{t \in (0, \infty)} \: t \: \parallel \Riem(g_t) 
\parallel_\infty \: < \: \infty$.

In three dimensions Perelman
has given important information about the long-time behavior
of Ricci flow solutions
\cite{Perelman (2002),Perelman (2003a),Perelman (2003b)},
which is especially relevant for topological purposes,
but the precise behavior is largely unknown. All known
compact three-dimensional Ricci flow solutions that exist
for all $t \in (0, \infty)$ are type-III, but it is not known whether
this is always the case. Hamilton had shown earlier that the
geometrization conjecture holds for such manifolds
\cite{Hamilton (1999)}.

This paper is concerned with the long-time behavior of 
$n$-dimensional Ricci flow
solutions, which we assume to be type-III.
Given a Ricci flow solution $g(\cdot)$ and a parameter $s > 0$,
there is another Ricci flow solution $g_s(\cdot)$ given by
$g_s(t) \: = \: s^{-1} \: g(st)$. The time interval $[a,b]$ for $g_s$
corresponds to the time interval $[sa, sb]$ for $g$.
Understanding the behavior of $g(t)$ for large $t$
amounts to understanding the behavior of
$g_s(\cdot)$ as $s \rightarrow \infty$. 

We show that in many examples
there is a limit as $s \rightarrow \infty$ of $g_s(\cdot)$, 
given by an expanding soliton $g_\infty(\cdot)$.
An expanding soliton has the scaling property that $g_\infty(t)$ differs 
from $t \: g_\infty(1)$ only by the action of a diffeomorphism $\eta_t$.
That the limit metric is expanding by the factor $t$
may seem contradictory to the fact that there are
compact Ricci
flow solutions that collapse, so we must explain in what
sense there is a limit and where it lives.

For concreteness, let us first discuss the
case of a locally homogeneous finite-volume $3$-manifold. The
lifted flow $\widetilde{g}(\cdot)$ on the universal cover $M$ has been extensively
studied; see Isenberg-Jackson
\cite{Isenberg-Jackson (1995)} and Knopf-McLeod
\cite{Knopf-McLeod (2001)}. In order to obtain
a limit $\widetilde{g}_\infty(\cdot)
\: = \: \lim_{s \rightarrow \infty} \widetilde{g}_s(\cdot)$ we use {\em pointed convergence}
of Ricci flows. 
Roughly speaking,  instead of comparing metrics on $M$ with respect to a fixed
coordinate system, we allow ourselves to transform
the metric $\widetilde{g}_s(t)$ by a $s$-dependent diffeomorphism. In effect,
we are choosing coordinates based on what an observer inside of the
manifold sees.

\begin{theorem} \label{intro1}
If $\widetilde{g}(\cdot)$ is the Ricci flow on a three-dimensional homogeneous
manifold that admits finite-volume quotients, and exists for all
$t \in (0, \infty)$, then there is a limit Ricci flow
$\widetilde{g}_\infty(\cdot) \: = \: \lim_{s \rightarrow \infty} \widetilde{g}_s(\cdot)$
which is an expanding soliton solution.
\end{theorem}

For each of the $3$-dimensional homogeneous classes there is
a unique limit soliton $\widetilde{g}_\infty(\cdot)$. It may be in a different
homogeneity class than the initial metric.  The expanding
solitons that we find are of type $\R^3$, $\R \times H^2$, $H^3$, $\Sol$ and $\Nil$.
If we start with an initial metric of type $\widetilde{\Isom^+(\R^2)}$ or 
$\widetilde{\SL(2, \R)}$ then we end up with an expanding soliton of type
$\R^3$ or $\R \times H^2$, respectively. In Section \ref{4D} we extend
Theorem \ref{intro1} to the four-dimensional homogeneous metrics considered by
Isenberg-Jackson-Lu \cite{Isenberg-Jackson-Lu (2005)}. Again we find
that there are limits $\widetilde{g}_\infty(\cdot) \: = \: \lim_{s \rightarrow \infty} \widetilde{g}_s(\cdot)$
given by expanding solitons.  

In these examples,
the metric $\widetilde{g}_\infty(t)$ gives $M$ the 
structure of a Riemannian submersion whose fiber is a
nilpotent Lie group and whose holonomy preserves the affine-flat structure
of the fiber.  The diffeomorphisms $\eta_t$ act fiberwise by means of
Lie group automorphisms.  This is related to the $\Nil$-structure described by
Cheeger-Fukaya-Gromov \cite{Cheeger-Fukaya-Gromov (1992)} for
collapse with bounded sectional curvature, and suggests that the
expanding solitons which are relevant for type-III solutions
may have a special structure.  Based on this, in Section
\ref{vectorbundles} we consider the expanding soliton equation
in the simplest case of a $\Nil$-structure,
namely when a manifold $M$ has a free isometric $\R^N$-action.

\begin{theorem} \label{intro2}
Let $M$ be the total space of a flat $\R^N$-vector bundle over a Riemannian
manifold $B$, with flat Riemannian metrics on the fibers.
Suppose that the fiberwise volume forms are preserved by
the flat connection.  Let $V(t)$ be the fiberwise
radial vector field $\frac{1}{2t} \sum_{i=1}^N x^i \: \frac{\partial}{\partial x_i}$.
Then the expanding soliton equation on $M$ becomes the equation for
a harmonic map $G \: : \: B \rightarrow  \SL(N, \R)/\SO(N)$ along with
the equation
\begin{equation} \label{introeqn}
R_{\alpha \beta} \: - \: \frac14 \: \Tr \left( G^{-1} \: G_{,\alpha} \: G^{-1} \: G_{,\beta}
\right) \: + \: \frac{1}{2t} \: g_{\alpha \beta} \: = \: 0
\end{equation}
on $B$.
\end{theorem}

In writing (\ref{introeqn}) we think of $G$ as taking value in 
positive-definite symmetric
$(N \times N)$-matrices.
Using this result, we give relevant examples of expanding solitons.

The results mentioned so far  mostly concern limit Ricci flows on
noncompact manifolds,  which may arise from Ricci flows on compact
manifolds upon taking the universal cover. 
One would also like to construct a limit flow for the
compact manifold.
Hamilton's compactness theorem gives sufficient conditions for a sequence
$\{g_k(\cdot)\}_{k=1}^\infty$ of pointed
 Ricci flow solutions to have a convergent subsequence
 \cite{Hamilton (1995)}. 
However, in order
to apply it one needs a uniform lower injectivity radius bound 
$ \inj_{g_k(t_0)} (p_k)  \: \ge \: i_0 \: > \:  0$.  In our case this precludes the
collapsing situation.  In order to obtain a limit flow in the collapsing case one
must consider Ricci flow on a larger class of spaces than smooth
manifolds.  One might try to
define Ricci flow on a Gromov-Hausdorff limit space, but this is not
very convenient.  Instead we will allow the limit Ricci flow to live on a space
which in a sense has the same dimension as the original manifold but which takes
the collapsing symmetry into account.  A convenient language is that
of Riemannian groupoids.  A Riemannian groupoid is an
\'etale groupoid equipped with an invariant Riemannian metric.
Riemannian groupoids have a history in
foliation theory, where they are used to describe
the transverse structure of Riemannian foliations; 
see Haefliger \cite{Haefliger (2001)} and references therein.  More
recently a similar notion was introduced by Petrunin and Tuschmann
under the name  ``megafold'' \cite[Appendix]{Petrunin-Tuschmann (1999)},
with application to collapsing in Riemannian geometry.
Two definitions were given in \cite[Appendix]{Petrunin-Tuschmann (1999)},
one in terms of topoi and one in terms of pseudogroups.  We prefer
the Riemannian groupoid language, but all three definitions are essentially
equivalent. We give an extension of Hamilton's compactness theorem to the 
case when
there is no positive lower bound on the injectivity radius. The limit Ricci flow will 
not be on a manifold but rather on
a groupoid.  

\begin{theorem} \label{intro3}
Let $\{(M_i, p_i, g_i(\cdot))\}_{i=1}^\infty$ be a sequence of
Ricci flow solutions on pointed $n$-dimensional
manifolds $(M_i, p_i)$. We assume that 
there are numbers $-\infty \: \le A \: < \: 0$ and $0 \: < \:\Omega
\: \le \: \infty$ so that \\
1. The Ricci flow solution $(M_i, p_i, g_i(\cdot))$ is defined on the
time interval $(A, \Omega)$. \\
2. For each $t \in (A,\Omega)$, $g_i(t)$ is a complete Riemannian metric
on $M_i$. \\
3. For each closed interval $I  \subset (A, \Omega)$
there is some $K_{I} \: < \: \infty$ so that $|\Riem(g_i)(x, t)| \: \le \:
K_{I}$ for all $x \in M_i$ and $t \in I$. \\

Then after passing to a subsequence,
the Ricci flow solutions $g_i (\cdot)$ converge
smoothly to 
a Ricci flow solution $g_\infty(\cdot)$ on a pointed 
$n$-dimensional \'etale groupoid 
$\left( G_\infty, O_{x_\infty} \right)$, defined again for 
$t \in (A, \Omega)$.
\end{theorem}

A result in this direction was proven by
Glickenstein \cite{Glickenstein (2003)} who constructed a limit flow
on a ball in a 
single tangent space; groupoids give a way of piecing these limits together
for various tangent spaces. Using the results of Section \ref{3d}, we show
that if $g(\cdot)$ is a Ricci flow on a
finite-volume locally homogeneous $3$-dimensional manifold, that
exists for all $t \in \infty$, then
$\lim_{s \rightarrow \infty} g_s(\cdot)$ exists and is an expanding
soliton on a $3$-dimensional \'etale groupoid.

Given $K > 0$, the space of pointed $n$-dimensional
Ricci flow solutions on manifolds with
$\sup_{t \in (0, \infty)} \: t \: \parallel \Riem(g_t) 
\parallel_\infty \: \le \: K$ is precompact among Ricci flows on
pointed $n$-dimensional
\'etale groupoids.  The closure ${\mathcal S}_{n,K}$ has an
$\R^+$-action that sends $g$ to $g_s$.  Understanding the
long-time behaviour of $n$-dimensional type-III Ricci flow solutions translates to 
understanding the  dynamics of the $\R^+$-action on
${\mathcal S}_{n,K}$, which seems to be an interesting problem.

The organization of the paper is as follows.  In Section \ref{deff} we
give some basic results about expanding solitons.  In Section
\ref{homogeneous} we consider the long-time behavior of Ricci
flow on homogeneous spaces of dimension one through four.
In Section \ref{vectorbundles} we look at the expanding soliton 
equation on a space with a free isometric $\R^N$-action and
reduce it to the harmonic-Einstein equations. In Section
\ref{rgroupoids} we recall basic facts about Riemannian
groupoids and give the extension of Hamilton's compactness
theorem. 

More detailed descriptions are at the beginnings of the sections.

I am grateful to Peng Lu for helpful discussions and for 
detailed explanations of his joint work in 
\cite{Isenberg-Jackson-Lu (2005)}.

\section{Expanding Solitons} \label{deff}

In this section we recall some basic properties of expanding solitons.
We also recall the definition of pointed convergence of a sequence of
Ricci flows. We define the rescaling $g_s(\cdot)$ of a Ricci flow solution
$g(\cdot)$ defined for $t \in (0, \infty)$. We show that if $\{g_s(\cdot)\}_{s > 0}$
has a limit as $s \rightarrow \infty$ then the limit $g_\infty(\cdot)$
is an expanding soliton.

\subsection{Definitions}

An {\em expanding soliton} on a manifold $M$
is a special type of Ricci flow solution
on a time interval $(t_0, \infty)$.
For convenience, we take $t_0 \: = \: 0$. Then the equation
for the time-dependent
Riemannian metric $g(t)$ and the time-dependent vector field
$V(t)$ is
\begin{equation} \label{expandingeqn1}
\Ric \: + \: \frac{{\mathcal L}_V g}{2} \: + \: \frac{g}{2t} \: = \: 0.
\end{equation}
Also, $V(t) \: = \: \frac{1}{t} \: V(1)$.
The corresponding Ricci flow is given by
\begin{equation} \label{expandingeqn2}
g(t) \: = \: t \: \eta_t^* g(1),
\end{equation}
where $\{\eta_t\}_{t>0}$ 
is the $1$-parameter family of diffeomorphisms generated 
by $\{V(t)\}_{t>0}$, normalized by $\eta_1 \: = \: \Id$. 
(If $M$ is noncompact then we assume that $V$ is such that we 
can solve for the $1$-parameter family.)
Conversely, given a solution to the time-independent equation
\begin{equation} \label{expandingeqn3}
\Ric \: + \: \frac{{\mathcal L}_V g}{2} \: + \: \frac{g}{2} \: = \: 0,
\end{equation}
put $V(t) \: = \: \frac{1}{t} \: V$, solve for
$\{\eta_t\}_{t>0}$ and put $g(t) \: = \: t \: \eta_t^* g$.
Then $(g(t), V(t))$ satisfies (\ref{expandingeqn1}).

If $(M_1, g_1(\cdot))$ and $(M_2, g_2(\cdot))$ are two expanding
soliton solutions, with associated diffeomorphisms
$\{\eta^{(1)}_t\}_{t>0}$ and $\{\eta^{(2)}_t\}_{t>0}$,
then the product flow
$(M_1 \times M_2, g_1(\cdot) + g_2(\cdot))$ is an expanding soliton
solution with $\eta_t \: = \: \left(
\eta^{(1)}_t, \eta^{(2)}_t \right)$.

Let $g(\cdot)$ be an expanding soliton on $M$. Suppose that $\Gamma$ is a
discrete group that acts on $M$
freely, properly discontinuously and isometrically
(with respect to the metrics $g(\cdot)$).
Then there is an quotient Ricci flow solution
$\overline{g}(\cdot)$ on $M/\Gamma$. 
If $\Gamma$ also preserves the vector
fields $V(\cdot)$ then $\overline{g}(\cdot)$ is an expanding
soliton, but this does not have to be the case.

\subsection{Expanding solitons as long-time limits} \label{def2}

Let $(M, p)$ be a connected manifold with a basepoint $p$.
Let $\{g(t)\}_{t \in (0, \infty)}$ be a Ricci flow solution on
$M$. We assume that 
for each $t > 0$, the pair $(M, g(t))$ is a complete Riemannian
manifold.
If $\{(M_i, p_i, g_i(\cdot))\}_{i=1}^\infty$ is a sequence of such
Ricci flow solutions
then there is a notion of {\em pointed convergence} to a limit
Ricci flow solution $(M_\infty, p_\infty, g_\infty(\cdot))$, as
considered in \cite{Hamilton (1995)}. In our case,
this means that one has \\
1. A
sequence of open subsets $\{U_j\}_{j=1}^\infty$ of $M_\infty$ containing
$p_\infty$, so that any compact subset of $M_\infty$ eventually lies in
all $U_j$, and \\
2. Time-independent diffeomorphisms $\phi_{i,j} \: : \: U_j 
\rightarrow V_{i,j}$ from $U_j$ to open subsets $V_{i,j} \subset M_i$, with
$\phi_{i,j} (p_\infty) \: = \: p_i$, so that\\
3. For all $j$,
$\lim_{i \rightarrow \infty} \phi_{i,j}^* g_i(\cdot) \: = \:
g_\infty (\cdot) \big|_{U_j}$ smoothly on $U_j \times [j^{-1}, j]$. \\

The compactness theorem of \cite{Hamilton (1995)} implies the
following. Suppose that \\
1. For each closed interval $I \subset (0, \infty)$ there is some
$K_I  < \infty$ 
so that $|\Riem|(x,t) \: \le \: K_I$ for all $x \in M_i$ and $t \in I$.\\
2. There are some $t_0 > 0$ and $i_0 > 0$ so that for all $i$, 
$\inj_{g_i(t_0)}(p_i) \: \ge \: i_0$. \\
Then $\{g_i(\cdot)\}_{i=1}^\infty$ has a convergent subsequence.

Given a $1$-parameter family $\{M, p, g_s(\cdot)\}_{s > 0}$ of Ricci flow
solutions, there is an analogous
notion of convergence as $s \rightarrow \infty$, i.e. for any sequence
$\{s_j\}_{j=1}^\infty$ converging to infinity
the sequence
$\{M, p, g_{s_j}(\cdot)\}_{j=1}^\infty$ converges 
and the limit is independent of the choice of $\{s_j\}_{j=1}^\infty$.
Hence it makes sense to talk about having
a limit solution
$\lim_{s \rightarrow \infty} (M, p, g_s(\cdot)) \: = \:
(M_\infty, p_\infty, g_\infty(\cdot))$

Now suppose that we have a type-III Ricci flow solution
$(M, p, g(\cdot)) $, meaning that
$\sup_{t \in (0, \infty)} \: t \: \parallel \Riem(g_t) \parallel_\infty 
\: < \: \infty$.
For any $s > 0$, there is a rescaled Ricci flow solution 
$(M, p, g_s(\cdot))$ given by $g_s(t) \: = \: s^{-1} \: g(st)$.
We will consider the convergence or subconvergence of 
$(M, p, g_s(\cdot))$ as $s \rightarrow \infty$. It is important to
note that although all of the Ricci flow solutions $(M, p, g_s(\cdot))$
live on the same manifold $M$, the notion of convergence is not that of 
smooth metrics on $M$. Instead,  we are 
interested in {\em pointed convergence} as defined above.

\begin{lemma}
If $\liminf_{t  \rightarrow \infty} t^{-\frac12} \: \inj_{g(t)}(p) > 0$ then any sequence
$\{s_i\}_{i=1}^\infty$ converging to infinity has a subsequence, which
we again denote by $\{s_i\}_{i=1}^\infty$, so that
$\lim_{i \rightarrow \infty} (M, p, g_{s_i}(\cdot)) \: = \: 
(M_\infty, p_\infty, g_{\infty}(\cdot))$ for some Ricci flow solution 
$(M_\infty, p_\infty, g_{\infty}(\cdot))$ defined for $t \in (0, \infty)$.
\end{lemma}
\begin{proof}
This is an immediate consequence of Hamilton's compactness theorem.
\end{proof}

We now consider what happens if there actually is a limit.

\begin{proposition} \label{limit}
If $\lim_{s \rightarrow \infty} (M, p, g_{s}(\cdot)) \: = \:
(M_\infty, p_\infty, g_\infty(\cdot))$ then
$(M_\infty, g_\infty(\cdot))$ is an expanding soliton.
\end{proposition}
\begin{proof}
Let ${\mathcal M}$ denote the space of pointed Riemannian metrics on 
$M_\infty$, with
the  topology of smooth convergence on compact subsets. 
The Ricci flow solution $g_\infty(\cdot)$ defines a smooth curve in ${\mathcal M}$.
Given $t, \alpha > 0$, we can formally write (modulo diffeomorphisms)
\begin{equation}
g_\infty(\alpha t) \: = \:
\lim_{s \rightarrow \infty} s^{-1} \: g(s \alpha t) \: = \:
\lim_{s \rightarrow \infty} \alpha \: s^{-1} \: g(s t) \: = \:
\alpha \: g_\infty(t).
\end{equation}
More precisely, for any $R > 0$ and $\epsilon > 0$ there is a pointed
(= basepoint-preserving) 
diffeomorphism
$\phi_{R,\epsilon}$ from the time-$\alpha t$ ball $B_R(p_\infty) \subset
M_\infty$ to a subset $V_{R, \epsilon} \subset M_\infty$ such that 
$\alpha \: \phi_{R,\epsilon}^* \: g_\infty(t) \big|_{V_{R, \epsilon}}$ is $\epsilon$-close 
in the smooth topology to 
$g_\infty(\alpha t) \big|_{B_R(p_\infty)}$. Taking the limit of
an appropriate sequence of the $\phi_{R,\epsilon}$'s, we obtain a 
pointed diffeomorphism $\phi \: : \: M_\infty \rightarrow M_\infty$ such that
$\alpha \: \phi^* g_\infty(t) \: = \: g_\infty(\alpha t)$.

Let $\Diff_{p_\infty}(M_\infty)$ denote the pointed
diffeomorphisms of $M_\infty$, again with the topology of smooth convergence
on compact subsets. We have shown that for all $t > 0$, the metric
$t^{-1} \: g_\infty(t)$ lies in the $\Diff_{p_\infty}(M_\infty)$-orbit of $g_\infty(1)$.
As in \cite{Bourguignon (1975)}, the $\Diff_{p_\infty}(M_\infty)$-orbit of $g_\infty(1)$
is the image of a proper embedding of the smooth infinite-dimensional manifold
$\Diff_{p_\infty}(M_\infty)/\Isom_{p_\infty}(g_\infty(1))$ in ${\mathcal M}$. 
(Strictly speaking the paper
\cite{Bourguignon (1975)} deals with compact manifolds.) Hence the
smooth curve $t \rightarrow t^{-1} \: g_\infty(t)$ defines a smooth curve in
$\Diff_{p_\infty}(M_\infty)/\Isom_{p_\infty}(g_\infty(1))$, which we can
lift to a smooth curve in $\Diff_{p_\infty}(M_\infty)$.

Thus we have found a smooth $1$-parameter family of pointed
diffeomorphisms $\{\eta_t\}_{t > 0}$ so that
(\ref{expandingeqn2}) is satisfied for $g_\infty(\cdot)$.
Letting  $\{V(t)\}_{t > 0}$ be the generator of $\{\eta_t\}_{t > 0}$,
equation (\ref{expandingeqn1}) is satisfied. Substituting
(\ref{expandingeqn2}) into (\ref{expandingeqn1}) gives
$t \: {\mathcal L}_{V(t)} g_\infty(1) \: = \: 
{\mathcal L}_{V(1)} g_\infty(1)$.
Hence we may assume that $V(t) \: = \: \frac{1}{t} \: V(1)$ and redefine
$\{\eta_t\}_{t > 0}$.
This proves the proposition.
\end{proof}

\section{Homogeneous solutions} \label{homogeneous}

In this section we consider the Ricci flow on homogeneous
Riemannian manifolds of dimension one through four that admit
finite-volume quotients and exist for all $t \in (0, \infty)$. 
In dimensions
one through three we show that in all cases there is a limit Ricci flow solution
$g_\infty(\cdot) \:  = \: \lim_{s \rightarrow \infty} \phi_s^* g_s(\cdot)$ given by
an expanding soliton. We compute the soliton metric explicitly.
In dimension four we show that this is also true for the homogeneous
metrics considered in \cite{Isenberg-Jackson-Lu (2005)}. The main
task in all of these cases is to construct appropriate diffeomorphisms
$\phi_s$.

A pointed Gromov-Hausdorff limit of a sequence of homogeneous
manifolds is still homogeneous \cite[Corollary on p. 66]{Gromov (1981)}.
Hence if $(M, p, g(\cdot))$ is a homogeneous Ricci flow solution then
assuming that the limit exists, we know that 
$(M_\infty, p_\infty, g_\infty(\cdot)) \: = \:
\lim_{s \rightarrow \infty} (M, p, g_s(\cdot))$ is also homogeneous.
However, the isometry group may change in the limit.

We now examine the long-time limits for homogeneous
Ricci flow solutions of dimensions one through four.
The manifolds that we consider are simply-connected
homogeneous spaces $M = G/H$, where $G$ is a transitive group
of diffeomorphisms 
of $M$ and $H$ is the isotropy subgroup, assumed to be compact. We will
assume that
$G$ is connected and 
unimodular, i.e. has a bi-invariant Haar measure.  This will be the
case if $M$ admits finite-volume quotients.  We take the basepoint $p$ to be
the identity coset $eH$. The Riemannian metrics that we consider on
$G/H$ will be left-invariant.

Given the manifold $M$, there are various groups
$G \subset \Diff(M)$ that act transitively on $M$ with compact
isotropy group. We will
take  {\em minimal} such groups, i.e. no proper subgroup of $G$ 
acts transitively on $M$. This allows for the widest class
of Ricci flow solutions. However, we must note that a compact
quotient of $M$ may be of the form $\Gamma \backslash M$ where $\Gamma$ is
a freely-acting discrete subgroup of some larger such group $G^\prime$
containing $G$. For this reason,
for the purposes of the
geometrization conjecture one generally takes $G$ to be a {\em maximal}
element among the groups of diffeomorphisms of $M$ that act transitively
with compact isotropy group \cite[\S 5]{Scott (1983)},
\cite[Chapter 3]{Thurston (1997)}.

Given a homogeneous Ricci flow solution $g(\cdot)$, the question is
whether we can find pointed diffeomorphisms $\{\phi_s\}_{s > 0}$ so that
there is a limit Ricci flow solution 
$g_\infty(\cdot) \: = \: \lim_{s \rightarrow \infty}
\phi_s^* \: g_s(\cdot)$, where $g_s(t) \: = \: \frac1s \: g(st)$.
By Proposition \ref{limit}, the limit will necessarily be an expanding
homogeneous soliton solution.

\begin{remark}
We will see examples of expanding solitons on Lie groups
$G$ with the property that the rescaling diffeomorphisms
$\{\eta_t\}_{t > 0}$ arise from a $1$-parameter group 
$\{a_t\}_{t>0}$ of automorphisms of $G$, by $\eta_t \: = \:
a_{t^{-1}}$. If so, let $\Gamma$ be a
discrete subgroup of $G$. Then $\Gamma \backslash G$ with the quotient 
metric $\overline{g}(t)$ is isometric to the result of
quotienting $(G, tg(1))$ on the left by the subgroup
$a_{t^{-1}}(\Gamma)$. Thus we can basically either think of the metric as
evolving, or of the discrete group as evolving. \\
\end{remark}

\subsection{One dimension} \label{1d}

The manifold $M$ is $\R$, with $(G, H) \: = \: (\R, \{e\})$. 
The basepoint is $0 \in \R$.
The Ricci flow solution $g(t)$ is constant in $t$, equaling a flat metric
$g_0$. Then $g_s(t) \: = \: s^{-1} \: g_0$.
Let $\phi_s$ be multiplication by $\sqrt{s}$ on $\R$. Then
$\phi_s^* g_s(t) \: = \: g_0$. Hence there is a limit as $s \rightarrow
\infty$ of $\phi_s^* g_s(\cdot)$ given by
$g_\infty(t) \: = \: g_0$. We note that this is an expanding soliton
solution, with $\eta_t$
being multiplication on $\R$ by $t^{- \frac12}$. 

The quotient $S^1 \: = \: \Z \backslash \R$ has the
constant Ricci flow solution $(S^1, \overline{g}(t))$.
We can consider $(S^1, \overline{g}(t))$ to be isometric to 
the quotient of $(\R, t g_0)$ by $a_{t^{-1}}(\Z)$, where
$a_t$ is the automorphism of $\R$ given by 
multiplication by $\sqrt{t}$.

\subsection{Two dimensions} \label{2d}

The possible homogeneous spaces are $S^2$, $\R^2$ and $H^2$. 
Their pairs $(G, H)$ are
$(\SO(3), \SO(2))$,
$(\R^2, \{e\})$ and 
$(\Isom^+(H^2), \SO(2))$.
The Ricci flow on $S^2$ has finite extinction time, so we
do not consider it further. The case $\R^2$ is a product case,
and so has already been covered.

For the $H^2$ case,
let $g_{0}$ be a complete constant-curvature metric on the plane
with $\Ric\left( g_{0}\right) =-cg_{0}$ for
some $c>0$. The Ricci flow solution
starting at $g_{0}$ is given by
$g\left( t\right) =\left( 1+2ct\right) g_{0}$.
Then 
$g_s(t) \: = \: s^{-1} \: \left( 1+2cst\right) g_{0}$.
Taking $\phi_s \: = \: \Id$,
there is a limit as $s \rightarrow \infty$ of
$\phi_s^* g_s(\cdot)$, given by
$g_\infty(t) \: = \: 2ct g_0$.
This is independent of $c$ and
is an expanding soliton solution with $V = 0$.

\subsection{Three dimensions} \label{3d}

A homogeneous
Ricci flow on $S^3$ or $S^2 \times \R$ has finite extinction time, so we
do not consider it further.
The homogeneous spaces $\R^3$ and $H^2 \times \R$ are
product cases. By the previous discussion, after appropriate rescaling
their Ricci flows have expanding soliton limits.

We now list the cases $M = G/H$ by the group $G$.

\subsubsection{$G = \Isom^+(H^3)$} \label{H3}

The group $G$ is the connected component of the identity in $\SO(3,1)$.
The subgroup $H$ is $\SO(3)$.
Let $g_{0}$ be a complete constant-curvature metric on $\R^3$ with 
$\Ric \left( g_{0}\right) =-cg_{0}$ for
some $c>0$. The Ricci flow solution starting at $g_0$ is given by
$g\left( t\right) =\left( 1+2ct\right) g_{0}$.
Then $g_s(t) \: = \: s^{-1} \: \left( 1+2cst\right) g_{0}$.
Taking $\phi_s \: = \: \Id$, there is a limit as $s \rightarrow \infty$
of  $\phi_s^* g_s(\cdot)$, given by
$g_\infty(t) \: = \: 2ct g_0$. This is independent of $c$ and
is an expanding soliton solution with $V = 0$. \\

The remaining cases have trivial isotropy group $H$, i.e. $M = G$.
It is known that $M$ admits a Milnor frame, i.e. a left-invariant
orthonormal frame field $\{X_1, X_2, X_3\}$ so that
$[X_i, X_j] \: = \: \sum_k c^k_{\: \: ij} X_k$ with
$c^k_{\: \: ij}$ vanishing unless $i$, $j$ and $k$ are mutually distinct.
In this basis, the nonzero components of the curvature tensor
are of the form $K_{ijij}$.
If $\{\theta^1, \theta^2, \theta^3\}$ is the dual orthonormal
coframe then the Ricci flow solution can be written in the form
\begin{equation}
g\left( t\right) =A\left( t\right) \left( \theta ^{1}\right) ^{2}+B\left(
t\right) \left( \theta ^{2}\right) ^{2}+C\left( t\right) \left( \theta
^{3}\right) ^{2}.
\end{equation}%
We write $A\left( 0\right)
=A_{0},B\left( 0\right) =B_{0}$ and $C\left( 0\right) =C_{0}.$

In what follows, we use computations from \cite{Isenberg-Jackson (1995)} and
\cite{Knopf-McLeod (2001)}. We note that the metrics in
\cite{Isenberg-Jackson (1995)} and \cite{Knopf-McLeod (2001)}
different by a constant. Our normalizations will be those of
\cite{Isenberg-Jackson (1995)}.
However, we will use a Milnor basis as in \cite{Knopf-McLeod (2001)}.
The simpler solutions are listed first.

\subsubsection{$G = \Sol$} \label{sol}

The group $G$ is a semidirect product $\R^2  \widetilde{\times} \R$, where
$\R$ acts on $\R^2$ by $z \cdot (x,y) \: = \: (e^z x, e^{-z} y)$. The
subgroup $H$ is trivial. After a change of basis, the
 Lie algebra relations are
$[X_2, X_3] \: = \: X_1$, $[X_3, X_1] \: = \: 0$ and 
$[X_1, X_2] \: = \: - \: X_3$. The $\R^2$-factor is spanned by
$X_1$ and $X_3$, and the $\R$-factor is spanned by $X_2$.

The metric is
\begin{equation}
g(t) \: = \: A(t) (\theta^1)^2 \: + \: B(t) (\theta^2)^2 \: + \: 
C(t) (\theta^3)^2,
\end{equation}
where
\begin{equation}
d\theta^1 \: = \: - \: \theta^2 \wedge \theta^3, \: \: \: \:
d\theta^2 \: = \: 0, \: \: \: \:
d\theta^3 \: = \: \theta^1 \wedge \theta^2.
\end{equation}
The sectional curvatures are 
\begin{eqnarray}
K_{12} &=&\frac{(A-C)^{2}-4C^{2}}{4ABC} \\
K_{23} &=& \frac{(A-C)^{2}-4A^{2}}{4ABC} \notag \\
K_{31} &=& \frac{(A+C)^2}{4ABC} \notag.
\end{eqnarray}%
The Ricci flow is given by
\begin{eqnarray}
\frac{dA}{dt} &=&\frac{C^{2}-A^{2}}{BC} \\
\frac{dB}{dt} &=&\frac{(A+C)^2}{AC} \notag \\
\frac{dC}{dt} &=&\frac{A^2 - C^2}{AB}. \notag
\end{eqnarray}%

From \cite{Knopf-McLeod (2001)} the large-$t$ asymptotics are
$\lim_{t \rightarrow \infty} A(t) \: = \: \lim_{t \rightarrow \infty} C(t)
\: = \: \sqrt{A_0 C_0}$ and $B(t) \sim 4t$.
Then 
\begin{equation}
g_s(t) \sim s^{-1} \: \sqrt{A_0 C_0}
\:  \left( (\theta^1)^2 \: + \: (\theta^3)^2 \right)
\: + \: 4t \: (\theta^2)^2.
\end{equation}

We take coordinates $(x,y,z)$ for $G$ 
in which $\theta^1 \: + \: \theta^3 \: = \:
e^{-z} \: dx$, $\theta^1 \: - \: \theta^3 \: = \:
e^{z} \: dy$ and $\theta^2 \: = \: dz$. Define diffeomorphisms
$\phi_s \: : \: \R^3 \rightarrow G$ by
\begin{equation}
\phi_s(x,y,z) \: = \: \left( (A_0 C_0)^{- \frac14} \: \sqrt{s} x, 
(A_0 C_0)^{- \frac14} \: \sqrt{s} y, z \right).
\end{equation}
Then there is a limit as $s \rightarrow \infty$ of $\phi_s^* g_s(\cdot),$ 
given by 
\begin{equation} \label{solmetric}
g_\infty(t) \: = \: (\theta^1)^2 \: + \: (\theta^3)^2
\: + \: 4t \: (\theta^2)^2.
\end{equation}
This is an expanding soliton solution with
$\eta_t(x,y,z) \: = \: \left( t^{-\frac12} x, t^{-\frac12} y, z \right)$.
Its geometry is a $\Sol$-geometry. We note that it is {\em not} a
gradient expanding soliton. The equation for the soliton also appeared in
\cite{Baird-Danielo (2005)}.

\begin{example}
An example of a $\Sol$-manifold is given by the total space of a $T^2$-bundle
over $S^1$ whose monodromy is a hyperbolic element of $\SL(2, \Z)$.
Geometrically, the long-time asymptotics of its Ricci flow amount to
shrinking the $T^2$ fiber by a factor of $\sqrt{t}$ and then
multiplying the overall metric by $t$.
\end{example}

\subsubsection{$\Nil$} \label{nil3}

The group $G$ is a nontrivial central $\R$-extension of
$\R^2$. The subgroup $H$ is trivial. The
Lie algebra relations are
$[X_2, X_3] \: = \: - \: X_1$ and $[X_3, X_1] \: = \:
[X_1, X_2] \: = \: 0$. The $\R$-factor is spanned by $X_1$, and
the $\R^2$-factor is spanned by $X_2$ and $X_3$.

The metric is
\begin{equation}
g(t) \: = \: A(t) (\theta^1)^2 \: + \: B(t) (\theta^2)^2 \: + \: 
C(t) (\theta^3)^2,
\end{equation}
where
\begin{equation}
d\theta^1 \: = \: \theta^2 \wedge \theta^3, \: \: \: \:
d\theta^2 \: = \: 0, \: \: \: \:
d\theta^3 \: = \: 0.
\end{equation}
The sectional curvatures are 
\begin{eqnarray}
K_{12} &=&\frac{A}{4BC} \\
K_{23} &=&- \: \frac{3A}{4BC} \notag \\
K_{31} &=& \frac{A}{4BC} \notag.
\end{eqnarray}%
The Ricci flow is given by
\begin{eqnarray}
\frac{dA}{dt} &=&- \: \frac{A^2}{BC} \\
\frac{dB}{dt} &=&\frac{A}{C} \notag \\
\frac{dC}{dt} &=&\frac{A}{B}. \notag
\end{eqnarray}%
The solution is
\begin{eqnarray}  \label{nilflow}
A &=&A_{0}\left( 1+\frac{3A_{0}}{B_{0}C_{0}}t\right) ^{-1/3} \\
B &=&B_{0}\left( 1+\frac{3A_{0}}{B_{0}C_{0}}t\right) ^{1/3} \notag \\
C &=&C_{0}\left( 1+\frac{3A_{0}}{B_{0}C_{0}}t\right) ^{1/3}. \notag
\end{eqnarray}
Then
\begin{eqnarray}
g_{s}\left( t\right) &=&s^{-1}A_{0}\left( 1+\frac{3A_{0}}{B_{0}C_{0}}%
s t\right) ^{-1/3}\left( \theta ^{1}\right) ^{2}+s^{-1}B_{0}\left( 1+%
\frac{3A_{0}}{B_{0}C_{0}}st\right) ^{1/3}\left( \theta ^{2}\right) ^{2} +
\\
&&s^{-1}C_{0}\left( 1+\frac{3A_{0}}{B_{0}C_{0}}st\right)
^{1/3}\left( \theta ^{3}\right) ^{2} \notag \\
&\sim &\left( \frac{A_{0}^{2}B_{0}C_{0}}{3}\right)
^{1/3} \: s^{-\frac43} \: t^{- \frac13} \:
\left( \theta ^{1}\right) ^{2} \: + \notag \\
&&\left( \frac{%
3A_{0}B_{0}^{2}}{C_{0}}\right)^{1/3} \: s^{-\frac23} \: t^{\frac13} \:
\left( \theta^{2}\right)^{2} \: + \: 
\left( \frac{%
3A_{0}C_{0}^{2}}{B_{0}}\right)^{1/3} \: s^{-\frac23} \: t^{\frac13} \:
\left(
\theta^{3}\right)^{2}.
\notag
\end{eqnarray}%

We take coordinates $(x,y,z)$ for $G$ 
in which $\theta^1 \: = \: dx \: + \: \frac12 \: y \: dz \: - \:
\frac12 \: z \: dy$,
$\theta^2 \: = \: dy$ and $\theta_3 \: = \: dz$.
Define diffeomorphisms
$\phi_s \: : \: \R^3 \rightarrow G$ by
\begin{equation}
\phi_s(x,y,z) \: = \: \left( \left( 9A_{0}^{2}B_{0}C_{0}\right)
^{-1/6}s^{2/3} x, 
\left( \frac{3A_{0}B_{0}^{2}}{C_{0}}\right)
^{-1/6}s^{1/3} y,
\left( \frac{3A_{0}C_{0}^{2}}{B_{0}}\right)
^{-1/6}s^{1/3} z \right).
\end{equation}
Then there is a limit as $s \rightarrow \infty$ of $\phi_s^* g_s(\cdot)$, 
given by 
\begin{equation} \label{nilmetric}
g_\infty(t) \: = \: 
\frac{1}{3t^{1/3}} \: \left( {\theta}^{1}\right)
^{2} \: + \: t^{1/3} \: \left( \left( {\theta}^{2}\right)^{2} \: + \:
\left( \theta^{3}\right) ^{2} \right).
\end{equation}
This is an expanding soliton solution with
$\eta_t(x,y,z) \: = \: \left( t^{-\frac23} x, t^{-\frac13}y, 
t^{-\frac13} z \right)$. Its geometry is a $\Nil$-geometry.
The equation for the soliton also appeared in
\cite{Baird-Danielo (2005)} and, implicitly, in
\cite{Lauret (2001)}.
 
\begin{example}
If $\Gamma$ is a lattice in $\Nil$, consider any locally homogeneous
Ricci flow $\overline{g}(\cdot)$ on 
$M = \Gamma \backslash \Nil$. As $t \rightarrow \infty$, 
$(M, \overline{g}(t))$ will
approach the quotient of $(\Nil, tg_\infty(1))$
by the subgroup
$a_{t^{-1}}(\Gamma^\prime)$, where $a_t$ is the automorphism of $\Nil$ given
by $a_t(x,y,z) \: = \: \left( t^{\frac23} x, t^{\frac13}y, 
t^{\frac13} z \right)$ and $\Gamma^\prime$ is a subgroup of $\Nil$ that
is isomorphic to $\Gamma$.
\end{example}

\subsubsection{ $G = \widetilde{\Isom^+(\R^2)}$} \label{isomr2}

The group $G$ is the universal cover of the orientation-preserving
isometries of $\R^2$. It is a semidirect product $\R^2 
\widetilde{\times} \R$, where
$\R$ acts on $\R^2$ by rotation. The subgroup $H$ is trivial.
The Lie algebra relations are
$[X_2, X_3] \: = \: X_1$, $[X_3, X_1] \: = \: X_2$ and $[X_1, X_2] \: = \: 0$. 
The $\R^2$-factor is spanned by $X_1$ and $X_2$, and the $\R$-factor
is spanned by $X_3$.

The compact quotients of $G$, as smooth manifolds, 
admit flat metrics. Because of this,
the group $G$ is generally not considered with regard to the
geometrization conjecture. Nevertheless, it is relevant for homogeneous
Ricci flow solutions.

The metric is
\begin{equation}
g(t) \: = \: A(t) (\theta^1)^2 \: + \: B(t) (\theta^2)^2 \: + \: 
C(t) (\theta^3)^2,
\end{equation}
where 
\begin{equation}
d\theta^1 \: = \: - \: \theta^2 \wedge \theta^3, \: \: \: \:
d\theta^2 \: = \: - \: \theta^3 \wedge \theta^1, \: \: \: \:
d\theta^3 \: = \: 0.
\end{equation}
The sectional curvatures are 
\begin{eqnarray}
K_{23} &=& \frac{(A+B)^{2}-4A^{2}}{4ABC} \\
K_{31} &=& \frac{(A+B)^{2}-4B^{2}}{4ABC} \notag \\
K_{12} &=&\frac{\left( A-B\right) ^{2}}{4ABC}. \notag
\end{eqnarray}%
The Ricci flow is given by
\begin{eqnarray}
\frac{dA}{dt} &=&-\frac{A^{2}-B^{2}}{BC} \\
\frac{dB}{dt} &=&-\frac{B^{2}-A^{2}}{AC} \notag \\
\frac{dC}{dt} &=&\frac{\left( A-B\right) ^{2}}{AB}.\notag
\end{eqnarray}%

From \cite{Knopf-McLeod (2001)}, there are limits
$\lim_{t \rightarrow \infty} A(t) \: = \: \lim_{t \rightarrow \infty} B(t)
\: = \: A_*$ and
$\lim_{t \rightarrow \infty} C(t) \: = \: C_*$, where
$A_* \: = \: \sqrt{A_0 B_0}$ and $C_* \: = \: 
\frac{C_0}{2} \left( \sqrt{\frac{A_0}{B_0}} \: + \: 
\sqrt{\frac{B_0}{A_0}} \right)$.
Then
\begin{equation}
g_{s}\left( t\right) \: \sim \: s^{-1} \: A_*  \left( \left( \theta
^{1}\right) ^{2}\: + \: 
\left( \theta ^{2}\right) \right)^{2} \: + \: 
s^{-1} \: C_* \left( \theta ^{3}\right) ^{2}.
\end{equation}

Define a diffeomorphism $\phi_s \: : \: \R^3 \rightarrow \: G$ by
\begin{equation}
\phi_s(x,y,z) \: = \: \alpha_s(x,y) \: \beta_s(z),
\end{equation}
where
$\alpha_s(x,y) \: = \: e^{\sqrt{s} (xX_1+yX_2)}$ and
$\beta_s(z) \: = \: e^{\sqrt{s} zX_3}$. Letting $h^{-1} dh$ denote the
Maurer-Cartan form on $G$, we have
\begin{align}
& \phi_s^* (h^{-1} dh) \: = \: \beta_s^{-1} \: \alpha_s^{-1} \: 
d\alpha_s \: \beta_s \: + \: \beta_s^{-1} \: d\beta_s \: = \:
\sqrt{s} \: \beta_s^{-1} \: (dx \: X_1 \: + \: dy \: X_2 ) \: \beta_s
\: + \: \sqrt{s} \: dz \: X_3 \: = \\
& \sqrt{s} \: \left( \cos(\sqrt{s} z) dx \: + \: \sin(\sqrt{s} z) dy \right)
\: X_1 \: + \:
\sqrt{s} \: \left(- \: \sin(\sqrt{s} z) dx \: + \: \cos(\sqrt{s} z) dy \right)
\: X_2 \: + \: \sqrt{s} \: dz \: X_3. \notag
\end{align}

If $\cdot_i$ denotes the $X_i$-component of an element of the Lie algebra 
then
\begin{equation}
\phi_s^* g_{s}\left( t\right) \: = \:
s^{-1} \: A(st) \: \left( \phi_s^* (h^{-1} dh) \right)_1^2 \: + \: 
s^{-1} \: B(st) \: \left( \phi_s^* (h^{-1} dh) \right)_2^2 \: + \: 
s^{-1} \: C(st) \: \left( \phi_s^* (h^{-1} dh) \right)_3^2.
\end{equation}
We see that there is a limit as $s \rightarrow \infty$ of
$\phi_s^* g_{s}\left( \cdot \right)$, given by
\begin{equation}
g_\infty(t) \: = \: A_* \: (dx^2 \: + \: dy^2) \: + \: C_* \: dz^2.
\end{equation}
This is a flat metric on $\R^3$ and, as we have seen, is an expanding
soliton solution.

\subsubsection{$\widetilde{\SL\left( 2,\mathbb{R}\right)}$} \label{sl2}

The group $G$ is the universal cover of $\SL(2, \R)$.
The subgroup $H$ is trivial.
The Lie algebra relations are
$[X_2, X_3] \: = \: - \:  X_1$, $[X_3, X_1] \: = \: X_2$ and 
$[X_1, X_2] \: = \: X_3$. 

The metric is
\begin{equation}
g(t) \: = \: A(t) (\theta^1)^2 \: + \: B(t) (\theta^2)^2 \: + \: 
C(t) (\theta^3)^2,
\end{equation}
where 
\begin{equation}
d\theta^1 \: = \: \theta^2 \wedge \theta^3, \: \: \: \:
d\theta^2 \: = \: - \: \theta^3 \wedge \theta^1, \: \: \: \:
d\theta^3 \: = \: - \: \theta^1 \wedge \theta^2.
\end{equation}
The sectional curvatures are 
\begin{eqnarray}
K_{23} &=& \frac{(B-C)^2 - A(3A+2B+2C)}{4ABC} \\
K_{31} &=& \frac{[A-(B-C)]^2 - 4B(B-C)}{4ABC} \notag \\
K_{12} &=&\frac{[A+(B-C)]^2+4C(B-C)}{4ABC}. \notag
\end{eqnarray}%
The Ricci flow is given by
\begin{eqnarray}
\frac{dA}{dt} &=&\frac{\left( B-C\right) ^{2}-A^{2}}{BC} \\
\frac{dB}{dt} &=&\frac{\left( C+A\right) ^{2}-B^{2}}{AC} \notag\\
\frac{dC}{dt} &=&\frac{\left( A+B\right) ^{2}-C^{2}}{AB}. \notag
\end{eqnarray}
From \cite{Knopf-McLeod (2001)}, the large-$t$ asymptotics are given by
$\lim_{t \rightarrow \infty} A(t) \: = \: A_* > 0$,
$B(t) \sim 2t$ and $C(t) \sim 2t$.
Then
\begin{equation} \label{asymp}
g_s(t) \: \sim \: s^{-1} \: A_* \: \left( \theta
^{1}\right) ^{2} \: + \: 2t \: \left( \left( \theta ^{2}\right)
^{2} \: + \: \left( \theta ^{3}\right) ^{2} \right).
\end{equation}

In order to find appropriate diffeomorphisms $\phi_s$ to
extract a limit, we use a bit of the geometry of
$\widetilde{\SL(2, \R)}$. Consider the $\R$-subgroup of 
$\widetilde{\SL(2, \R)}$ generated by $X_1$. For any $t > 0$, 
the metric on $\widetilde{\SL(2, \R)}$ given by (\ref{asymp})
is the total space of a Riemannian submersion
$\pi \: : \: \widetilde{\SL(2, \R)} \rightarrow
\R \backslash \widetilde{\SL(2, \R)}$, with 
$\R \backslash \widetilde{\SL(2, \R)}$ having
a metric of constant negative curvature. Furthermore, the fibers of
the submersion are totally geodesic lines. The idea is to
construct something like Fermi coordinates around the
fiber through the identity element $e \in \widetilde{\SL(2, \R)}$.

\begin{lemma}
Any element $h \in \widetilde{\SL(2, \R)}$ can be written uniquely
as $h \: = \: e^{aX_2 + bX_3} \: e^{cX_1}$ for some
$a,b,c \in \R^3$.
\end{lemma}
\begin{proof}
Choosing the metric
$\left( \theta
^{1}\right) ^{2} \: + \: \left( \theta ^{2}\right)
^{2} \: + \: \left( \theta ^{3}\right) ^{2}$ on 
$\widetilde{\SL(2, \R)}$ for concreteness, one can check that
for any $a,b \in \R$,  the
curve $v \rightarrow e^{avX_2 + bvX_3}$ is a geodesic in
$\widetilde{\SL(2, \R)}$. Clearly it is horizontal at
the identity element $e \in \widetilde{\SL(2, \R)}$. Hence it is horizontal
for all $v$ and represents the horizontal lift of a
geodesic in $\R \backslash \widetilde{\SL(2, \R)}$.
Now given an element $h \in \widetilde{\SL(2, \R)}$, consider the
unique geodesic
$\gamma \: : \: [0,1] \rightarrow \R \backslash \widetilde{\SL(2, \R)}$
with $\gamma(0) \: = \: \pi(e)$ and $\gamma(1) \: = \: \pi(h)$. 
Lift it to a horizontal 
geodesic $\widehat{\gamma}(v) \: = \: e^{avX_2 + bvX_3}$. As 
$\pi(\widehat{\gamma}(1)) \: = \: \pi(h)$, there is a unique $c \in \R$ so that
$h \: = \: e^{cX_1} \: \widehat{\gamma}(1)$.

This shows that $h$ can be written as uniquely as
$e^{cX_1} \: e^{aX_2 + bX_3}$ for some $a,b,c \in \R$. As
\begin{equation}
e^{cX_1} \: e^{aX_2 + bX_3} \: = \: 
e^{(\cos(c)a - \sin(c) b)X_2 + (\sin(c)a + \cos(c) b)X_3} \: e^{cX_1},
\end{equation}
the lemma follows.
\end{proof}

Define a diffeomorphism $\phi_s \: : \: \R^3 \rightarrow \: G$ by
\begin{equation}
\phi_s(x,y,z) \: = \: \alpha(y,z) \: \beta_s(x),
\end{equation}
where
$\alpha(y,z) \: = \: e^{yX_2+zX_3}$ and
$\beta_s(x) \: = \: e^{\sqrt{\frac{s}{A_*}} \: xX_1}$. 
Letting $h^{-1} dh$ denote the
Maurer-Cartan form on $G$, we have
\begin{eqnarray}
 \phi_s^* (h^{-1} dh) \: & = & \: \beta_s^{-1} \: \alpha^{-1} \: 
d\alpha \: \beta_s \: + \: \beta_s^{-1} \: d\beta_s \\
& = & \:\beta_s^{-1} \: \alpha^{-1} \: 
d\alpha \: \beta_s
\: + \: \sqrt{\frac{s}{A_*}} \: dx \: X_1. \notag
\end{eqnarray}

If $\cdot_i$ denotes the $X_i$-component of an element of the Lie algebra then
\begin{align}
\phi_s^* g_{s}\left( t\right) \: = \:
& s^{-1} \: A(st) \: \left( \phi_s^* (h^{-1} dh) \right)_1^2 \: + \: 
s^{-1} \: B(st) \: \left( \phi_s^* (h^{-1} dh) \right)_2^2 \: + \: 
s^{-1} \: C(st) \: \left( \phi_s^* (h^{-1} dh) \right)_3^2 \\
 = \: &  
s^{-1} \: A(st) \: \left( (\alpha^{-1} d\alpha)_1 \: + \:
\sqrt{\frac{s}{A_*}} \: dx
\right)^2 \: + \notag \\
& s^{-1} \: B(st) \: 
\left( \beta_s^{-1} \: \alpha^{-1} \: d\alpha \: \beta_s \right)_2^2 \: + \: 
s^{-1} \: C(st) \: \left( \beta_s^{-1} \: \alpha^{-1} \: d\alpha \: \beta_s
\right)_3^2. \notag
\end{align}
As conjugation by $\beta_s$ amounts to a rotation in the $(y,z)$-plane,
we see that there is a limit as $s \rightarrow \infty$ of
$\phi_s^* g_{s}\left( \cdot \right)$, given by
\begin{equation}
g_\infty(t) \: = \: dx^2 \: + 2t \: \left(
\left( \alpha^{-1} d\alpha \right)_2^2 \: + \: 
\left( \alpha^{-1} d\alpha \right)_3^2 \right).
\end{equation}
The pullback under $\pi \circ \alpha$ of the metric on 
$\R \backslash \widetilde{\SL(2, \R)}$ 
is the same as $\left( \alpha^{-1} d\alpha \right)_2^2 \: + \: 
\left( \alpha^{-1} d\alpha \right)_3^2$.
Hence $g_\infty(\cdot)$
is the expanding soliton on $\R \times H^2$, with
$\eta_t(x,y,z) \: = \: \left( t^{-\frac12} x, y, z \right)$.

\begin{example} \label{sl2example}
An example of a locally homogeneous $\widetilde{\SL(2, \R)}$-geometry
is the unit sphere bundle $S^1 \Sigma$ of a closed hyperbolic surface $\Sigma$.
Let $\overline{g}(\cdot)$ be its Ricci flow. In the most direct picture,
the manifolds $\left( S^1 \Sigma, t^{-1} \overline{g}(\cdot) \right)$ have
a Gromov-Hausdorff limit, as $t \rightarrow \infty$, given by
the rescaling of $\Sigma$ which has constant sectional curvature
$- \: \frac12$. For a more refined picture,
let $\pi \: : \: S^1 \Sigma \rightarrow \Sigma$ be the projection map. 
Given $p^\prime \in \Sigma$, let $B$ be a small ball around
$p^\prime$ in $\Sigma$. Then $\pi^{-1} (B)$ is diffeomorphic to
$S^1 \times B$. Consider the restriction of $t^{-1} \: \overline{g}(t)$
to $\pi^{-1} (B)$ 
and then its pullback to the universal cover $\widetilde{\pi^{-1} (B)}$.
Taking a basepoint $p \in \widetilde{\pi^{-1} (B)}$ over
$p^\prime$, the pointed limit as $t \rightarrow \infty$ of the
metric on $\widetilde{\pi^{-1} (B)}$ will be isometric to
$\R$ times a ball of constant sectional curvature
$- \: \frac12$. In effect, as the circle fibers of the manifold
$\left( S^1 \Sigma, t^{-1} \: \overline{g}(t) \right)$ shrink, the
local geometry becomes more and more product-like.
\end{example}

\subsection{Four dimensions} \label{4D}

We first list the four-dimensional homogeneous spaces $G/H$ with
{\em maximal} groups $G$ (acting transitively with compact isotropy
group) that admit finite-volume quotients
\cite{Wall (1986)}.
Besides product cases, they are
\begin{eqnarray} \label{wall}
\underline{G} & \underline{H} & \underline{G/H} \\ 
\SO(5) & \SO(4) & S^4  \notag \\
\SU(3) & U(2) & \C P^2  \notag \\
\Isom^+(H^4) & \SO(4) & H^4  \notag \\
\SU(2, 1) & U(2) & \C H^2  \notag \\
\R^2 \widetilde{\times} \SL(2, \R) & \SO(2) & F_4 \notag \\
(\C \times \R) \widetilde{\times} \C^* & \SO(2) & \Sol_0 \notag \\
\Nil^4 & \{e\} & \Nil^4 \notag \\
\Sol^4_{m,n} & \{e\} & \Sol^4_{m,n} \notag \\
\Sol^4_1 & \{e\} & \Sol^4_1. \notag
\end{eqnarray}

In the definition of $\Sol_0$, the action of $\C^*$ on
$\C \times \R$ is given by $\lambda \cdot (a,b) \: = \:
(\lambda a, |\lambda|^{-2} b)$. The group $\Nil^4$ is the
semidirect product $\R^3 \widetilde{\times} \R$, where
$\R$ acts on $\R^3$ by
$\beta(r) \: = \:
\begin{pmatrix}
1 & r & \frac{r^2}{2} \\
0 & 1 & r \\
0 & 0 & 1
\end{pmatrix}$. The group $\Sol^4_{m,n}$ is the
semidirect product $\R^3 \widetilde{\times} \R$, where
$\R$ acts on $\R^3$ by
$r \cdot (x,y,z) \: = \: 
\left( e^{ar}x, e^{br}y, e^{cr}z \right)$.
Here $a > b > c$ are real, $a+b+c=0$ and
$e^a$, $e^b$, $e^c$ are the roots of
$\lambda^3 - m \lambda^2 + n \lambda - 1 = 0$ with
$m,n \in \Z^+$. Finally, 
$\Sol^4_1 \: = \:
\left\{ 
\begin{pmatrix}
1 & b & c \\
0 & \alpha & a \\
0 & 0 & 1
\end{pmatrix}
\: : \: \alpha, a, b, c \in \R, \: \alpha > 0
\right\}$.

There is also a list of nonmaximal geometries 
\cite[Theorem 3.1]{Wall (1986)} but we do not consider it here.

A homogeneous Ricci flow on $S^4$ or $\C P ^2$ has finite
extinction time, so we do not consider it further. 

\subsubsection{$G = \Isom^+(H^4)$} \label{h4}

The group $G$ is the connected component of the identity in $\SO(4,1)$.
The subgroup $H$ is $\SO(4)$.
Let $g_{0}$ be a complete constant-curvature metric on $\R^4$ with 
$\Ric \left( g_{0}\right) =-cg_{0}$ for
some $c>0$. The Ricci flow solution starting at $g_0$ is given by
$g\left( t\right) =\left( 1+2ct\right) g_{0}$.
Then $g_s(t) \: = \: s^{-1} \: \left( 1+2cst\right) g_{0}$.
Taking $\phi_s \: = \: \Id$, there is a limit as $s \rightarrow \infty$
of  $\phi_s^* g_s(\cdot)$, given by
$g_\infty(t) \: = \: 2ct g_0$. This is independent of $c$ and
is an expanding soliton solution with $V = 0$.

\subsubsection{$G = \SU(2,1)$} \label{ch2}

The subgroup $H$ is $U(2)$.
Let $g_{0}$ be a complete metric on $\C^2$ with constant
holomorphic sectional curvature and with 
$\Ric \left( g_{0}\right) =-cg_{0}$ for
some $c>0$. The Ricci flow solution starting at $g_0$ is given by
$g\left( t\right) =\left( 1+2ct\right) g_{0}$.
Then $g_s(t) \: = \: s^{-1} \: \left( 1+2cst\right) g_{0}$.
Taking $\phi_s \: = \: \Id$, there is a limit as $s \rightarrow \infty$
of  $\phi_s^* g_s(\cdot)$, given by
$g_\infty(t) \: = \: 2ct g_0$. This is independent of $c$ and
is an expanding soliton solution with $V = 0$.

\subsubsection{$G = \R^2 \widetilde{\times} \SL(2, \R)$} \label{f4}

The quotient space $M = G/H$ fibers homogeneously over $H^2$, with fiber
$\R^2$. Any left-invariant metric defines a homogeneous
Riemannian submersion $M \rightarrow H^2$.
The isotropy group $\SO(2)$ acts isometrically on the Riemannian submersion,
rotating the $\R^2$-fiber containing the basepoint $\star$. By the
rotational symmetry, the
curvature tensor of the Riemannian submersion, a horizontal
$2$-form with values in the vertical tangent bundle, must vanish
at $*$. Then by homogeneity, it must vanish everywhere.  Thus
the horizontal space is integrable.  It follows that the Riemannian
submersion is of the type considered in Section \ref{vectorbundles},
so we can apply the results of that section; see Example
\ref{F4}. In particular, the
homogeneous metric on $M$ is specified by the relative size of the
fiberwise metric $g_{ij}$ on $\R^2$ and the relative size of the base
metric $g_{\alpha \beta}$ on $H^2$.
The Ricci curvature calculation of (\ref{ricci2}) shows that under the
Ricci flow, $g_{ij}$ is constant in $t$ and $g_{\alpha \beta}$ 
increases linearly in $t$. Taking $\phi_s$ to be multiplication by
$\sqrt{s}$ in the $\R^2$-fibers, there is a limit
as $s \rightarrow \infty$ of  $\phi_s^* g_s(\cdot)$, given by
the expanding soliton solution on F4. It satisfies the
harmonic-Einstein equations of Proposition \ref{harmoniceinstein}. \\

The remaining cases can be seen as Ricci flows on certain unimodular Lie
groups.  For example, the case $\Sol_0$ can be viewed as Ricci flow on the
Lie group
$\R^3 \widetilde{\times} \R$, where
$\R$ acts on $\R^3$ by
$\delta(r) \: = \:
\begin{pmatrix}
e^r & 0 & 0 \\
0 & e^r & 0 \\
0 & 0 & e^{-2r}
\end{pmatrix}$.
In \cite{Isenberg-Jackson-Lu (2005)} the Ricci flow was considered
on a class of metrics on four-dimensional unimodular Lie groups
that have the
property that their Ricci flow ``diagonalizes''.  The groups are listed
as A1-A10 in \cite{Isenberg-Jackson-Lu (2005)}.
They include some that do not admit
finite-volume quotients. In what follows we will use the calculations of
\cite{Isenberg-Jackson-Lu (2005)}. We now 
go through the cases A1-A10 in order.

\subsubsection{A1} \label{A1} This is flat $\R^4$.

\subsubsection{A2} \label{A2} The group $G$ is a semidirect product
$\R^3 \widetilde{\times} \R$, where $\R$ acts on $\R^3$ by
$r \cdot(x,y,z) \: =\: \left( e^r x, e^{kr} y, e^{-(k+1)r}z \right)$.
Here $k$ is a free parameter. Special case are
$\Sol^4_0$ and $\Sol^4_{m,n}$. The nonzero Lie algebra relations are
$[X_1, X_4] = X_1$, $[X_2, X_4] = kX_2$ and $[X_3, X_4] = 
-(k+1) X_3$. The $\R$-factor is spanned by $X_4$.

The metric is
\begin{equation} \label{metric}
g(t) \: = \: A(t) (\theta^1)^2 \: + \: B(t) (\theta^2)^2 \: + \: 
C(t) (\theta^3)^2 \: + \: D(t) (\theta^4)^2,
\end{equation} 
where
\begin{equation}
d\theta^1 \: = \: - \: \theta^1 \wedge \theta^4, \: \: \: \:
d\theta^2 \: = \: - \: k \: \theta^2 \wedge \theta^4, \: \: \: \:
d\theta^3 \: = \: (k+1) \: \theta^3 \wedge \theta^4, \: \: \: \:
d\theta^4 \: = \: 0.
\end{equation}
We note that here, and in the cases that follow, the metric
(\ref{metric}) is not the
most general homogeneous metric on $G/H$. However, it is a metric
for which the Ricci flow ``diagonalizes''.

The Ricci flow is given by
\begin{align}
A(t) \: & = \: A_0, \\
B(t) \: & = \: B_0, \notag \\ 
C(t) \: & = \: C_0, \notag \\
D(t) \: & = \: D_0 \: + \: 4(k^2+k+1)t. \notag
\end{align}
Then 
\begin{equation}
g_s(t) \: \sim \: 
s^{-1} \: A_0 (\theta^1)^2 \: + \: 
s^{-1} \: B_0 (\theta^2)^2 \: + \: 
s^{-1} \: C_0 (\theta^3)^2 \: + \: 4(k^2+k+1)t \: (\theta^4)^2.
\end{equation}
Define diffeomorphisms $\phi_s \: : \: \R^4 \rightarrow G$ by
$\phi_s(x,y,z,r) \: = \: \left(
A_0^{- \frac12} \: \sqrt{s} x,
B_0^{- \frac12} \: \sqrt{s} y,
C_0^{- \frac12} \: \sqrt{s} z,
r \right)$.
Then there is a limit as $s \rightarrow \infty$ of
$\phi_s^* g_s(\cdot)$ given by
\begin{equation}
g_\infty(t) \: = \: 
(\theta^1)^2 \: + \: 
(\theta^2)^2 \: + \: 
(\theta^1)^2 \: + \: 4(k^2+k+1)t \: (\theta^4)^2.
\end{equation}
This is an expanding soliton solution with
$\eta_t(x,y,z,r) \: = \:
\left( t^{-\frac12} x, t^{-\frac12} y, t^{-\frac12} z, r \right)$.

\subsubsection{A3} \label{A3} The group $G$ is a semidirect product
$\R^3 \widetilde{\times} \R$, where $\R$ acts on $\R^3$ by
$\epsilon(r) \: = \: 
\begin{pmatrix}
e^{kr} \cos(r) & e^{kr} \sin(r) & 0 \\
 - \: e^{kr} \sin(r) &  e^{kr} \cos(r) & 0 \\
0 & 0 &  e^{-2kr}
\end{pmatrix}$. The nonzero Lie algebra relations are
$[X_1, X_4] \: = \: k X_1 \: + \: X_2$, 
$[X_2, X_4] \: = \: - \: X_1 \: + \: k X_2$ and
$[X_3, X_4] \: = \: -2k X_3$.
Here $k$ is a nonzero number. The $\R$-factor is spanned by $X_4$.

The metric is
\begin{equation}
g(t) \: = \: A(t) (\theta^1)^2 \: + \: B(t) (\theta^2)^2 \: + \: 
C(t) (\theta^3)^2 \: + \: D(t) (\theta^4)^2,
\end{equation} 
where
\begin{equation}
d\theta^1 \: = \: - \: k \theta^1 \wedge \theta^4 \: + \:
\theta^2 \wedge \theta^4, \: \: \: \:
d\theta^2 \: = \: - \: \theta^1 \wedge \theta^4 \: - \:
k \theta^2 \wedge \theta^4, \: \: \: \:
d\theta^3 \: = \: 2k \: \theta^3 \wedge \theta^4, \: \: \: \:
d\theta^4 \: = \: 0.
\end{equation}
The Ricci flow asymptotics are
\begin{align}
A(t) \: & \sim \: \sqrt{A_0 B_0}, \\
B(t) \: & \sim \: \sqrt{A_0 B_0}, \notag \\ 
C(t) \: & = \: C_0 \notag \\
D(t) \: & \sim \: 12 k^2 t. \notag
\end{align}
Then 
\begin{equation}
g_s(t) \: \sim \: 
s^{-1} \: \sqrt{A_0 B_0} \left(
(\theta^1)^2 \: + \: 
(\theta^2)^2 \right) \: + \: 
s^{-1} \: C_0 \: (\theta^3)^2 \: + \: 12 k^2 t \: (\theta^4)^2.
\end{equation}

Define a diffeomorphism $\phi_s \: : \: \R^4 \rightarrow G$ by
\begin{equation}
\phi_s(x,y,z,r) \: =\: \alpha_s(x,y,z) \: \beta(r), 
\end{equation}
where
$\alpha_s(x,y,z) \: = \: e^{(A_0 B_0)^{- \frac14} \sqrt{s} 
(xX_1+yX_2) + C_0^{-\frac12} \sqrt{s} z X_3}$ and 
$\beta(r) \: = \: e^{k^{-1}rX_4}$.
Letting $h^{-1} dh$ denote the
Maurer-Cartan form on $G$, we have
\begin{align}
& \phi_s^* (h^{-1} dh) \: = \: \beta^{-1} \: \alpha_s^{-1} \: 
d\alpha_s \: \beta \: + \: \beta^{-1} \: d\beta \: = \\
& (A_0 B_0)^{- \frac14} \: \sqrt{s} \:
\beta^{-1} \: (dx \: X_1 \: + \: dy \: X_2) \: \beta \: + \: 
C_0^{-\frac12} \sqrt{s} \: dz \: \beta^{-1} \: X_3 \: \beta
\: + \: k^{-1} \: dr \: X_4 \: = \notag \\
&  (A_0 B_0)^{- \frac14} \:
\sqrt{s} \: e^{r} \: \left( \cos(k^{-1}
\sqrt{s} z) dx \: - \: \sin(k^{-1} \sqrt{s} z) dy \right)
\: X_1 \: + \notag \\
& (A_0 B_0)^{- \frac14} \:
\sqrt{s} \: e^{r} \: \left(\sin(k^{-1}
\sqrt{s} z) dx \: + \: \cos(k^{-1} \sqrt{s} z) dy \right)
\: X_2 \: + \notag \\
& C_0^{-\frac12} \: \sqrt{s} \: e^{- 2r} \:
dz \: X_3 \: + \: k^{-1} \: dr \: X_4. \notag
\end{align}

If $\cdot_i$ denotes the $X_i$-component of an element of the Lie algebra 
then
\begin{align}
\phi_s^* g_{s}\left( t\right) \: = \:
& s^{-1} \: A(st) \: \left( \phi_s^* (h^{-1} dh) \right)_1^2 \: + \: 
s^{-1} \: B(st) \: \left( \phi_s^* (h^{-1} dh) \right)_2^2 \: + \: 
s^{-1} \: C(st) \: \left( \phi_s^* (h^{-1} dh) \right)_3^2 \: + \\
& s^{-1} \: D(st) \: \left( \phi_s^* (h^{-1} dh) \right)_4^2. \notag
\end{align}
We see that there is a limit as $s \rightarrow \infty$ of
$\phi_s^* g_{s}\left( \cdot \right)$, given by
\begin{equation}
g_\infty(t) \: = \: e^{2r} \: (dx^2 \: + \: dy^2) \: + \: 
e^{-4r} \: dz^2 \: + \: 12t \: dr^2.
\end{equation}
This is an expanding soliton with
$\eta_t(x,y,z,r) \: = \: \left( t^{-\frac12} x, t^{-\frac12} y,
t^{-\frac12} z, r \right)$. It has $\Sol^4_0$-symmetry.

\subsubsection{A4} \label{A4}

This is a product case $G \: = \: \Nil^3 \times \R$.

\subsubsection{A5} \label{A5}
The group $G$ is a semidirect product
$\R^3 \widetilde{\times} \R$, where $\R$ acts on $\R^3$ by
$\phi(r) \: = \: 
\begin{pmatrix}
e^{-\frac{r}{2}} &  r e^{-\frac{r}{2}}  & 0 \\
0 &  e^{-\frac{r}{2}} & 0 \\
0 & 0 &  e^r
\end{pmatrix}$. The nonzero Lie algebra relations are
$[X_1, X_4] \: = \: - \: \frac12 \: X_1 \: + \: X_2$, 
$[X_2, X_4] \: = \: - \: \frac12 \:  X_2$ and
$[X_3, X_4] \: = \: X_3$. The $\R$-factor is spanned by $X_4$.

The metric is
\begin{equation}
g(t) \: = \: A(t) (\theta^1)^2 \: + \: B(t) (\theta^2)^2 \: + \: 
C(t) (\theta^3)^2 \: + \: D(t) (\theta^4)^2,
\end{equation} 
where
\begin{equation}
d\theta^1 \: = \: \frac12 \: \theta^1 \wedge \theta^4, \: \: \: \:
d\theta^2 \: = \: - \: \theta^1 \wedge \theta^4 \: + \:
\frac12 \: \theta^2 \wedge \theta^4, \: \: \: \:
d\theta^3 \: = \: - \: \theta^3 \wedge \theta^4, \: \: \: \:
d\theta^4 \: = \: 0.
\end{equation}
The Ricci flow asymptotics are
\begin{align}
A(t) \: & \sim \: 2 \lambda (\ln{t})^{\frac12}, \\
B(t) \: & \sim \: 3 \lambda (\ln{t})^{-\frac12}, \notag \\ 
C(t) \: & = \: C_0 \notag \\
D(t) \: & \sim \: 3 t. \notag
\end{align}
Then 
\begin{equation}
g_s(t) \: \sim \: 
2 \lambda s^{-1} \: (\ln(st))^{\frac12} \:
(\theta^1)^2  \: + \: 3 \lambda s^{-1} \: (\ln(st))^{-\frac12}
(\theta^2)^2 \: + \: 
s^{-1} \: C_0 \: (\theta^3)^2 \: + \: 3 t \: (\theta^4)^2.
\end{equation}

We take coordinates $(x,y,z,r)$ for $G$ in which
$\theta^1 \: = \: dx \: + \: \frac12 \: x dr$, $\theta^2 \: = \: dy \: + \:
\left( \frac{y}{2} - x \right) dr$, 
$\theta^3 \: = \: dz \: - \:  z dr$ and $\theta^4 \: = \: dr$.
Define $\phi_s \: : \: \R^4 \rightarrow \R^4$ by
\begin{equation}
\phi_s(x, y, z, r) \: = \: 
\left( \left( \frac{s}{2\lambda(\ln s)^{\frac12}} \right)^{\frac12} \: x,
\left( \frac{s(\ln s)^{\frac12}}{3\lambda} \right)^{\frac12} \: y,
\left( \frac{s}{C_0} \right)^{\frac12} \: z, \: r \right).
\end{equation}
Then $\lim_{s \rightarrow \infty} s^{-1} \phi_s^* g(st) \: = \: 
g_\infty(t)$, where
\begin{equation}
g_\infty(t) \: = \: (dx \: + \: \frac12 \: x \: dr)^2 \: + \:
(dy \: + \: \frac12 \: y \: dr)^2 \: + \: (dz \: - \:
z dr)^2 \: + \:
3t \: dr^2.
\end{equation}
This is an expanding soliton with
$\eta_t(x,y,z,r) \: = \: \left( t^{-\frac12} x, t^{-\frac12} y,
t^{-\frac12} z, r \right)$. It has $\Sol^4_0$-symmetry.

\subsubsection{A6}  \label{A6}

The group $G$ is the four-dimensional nilpotent Lie group whose
nonzero Lie algebra relations are
$[X_1, X_4] \: = \: X_2$ and $[X_2, X_4] \: = \: X_3$.

The metric is
\begin{equation}
g(t) \: = \: A(t) (\theta^1)^2 \: + \: B(t) (\theta^2)^2 \: + \: 
C(t) (\theta^3)^2 \: + \: D(t) (\theta^4)^2,
\end{equation} 
where
\begin{equation}
d\theta^1 \: = \: 0, \: \: \: \:
d\theta^2 \: = \: - \: \theta^1 \wedge \theta^4, \: \: \: \:
d\theta^3 \: = \: - \: \theta^2 \wedge \theta^4, \: \: \: \:
d\theta^4 \: = \: 0.
\end{equation}
Put $E_0 \: = \: \frac{B_0}{A_0 D_0}$ and
$F_0 \: = \: \frac{C_0}{B_0 D_0}$.
The Ricci flow solution is
\begin{align}
A(t) \: & = \: A_0 \: (3E_0t \: + \: 1)^{\frac13}, \\
B(t) \: & = \: B_0 \: (3E_0t \: + \: 1)^{-\frac13} \:
(3F_0t \: + \: 1)^{\frac13}, \notag \\ 
C(t) \: & = \: C_0 \: (3F_0t \: + \: 1)^{- \frac13} \notag \\
D(t) \: & = \: D_0 \: (3E_0t \: + \: 1)^{\frac13} \: 
(3F_0t \: + \: 1)^{\frac13}. \notag
\end{align}
Then 
\begin{align}
g_s(t) \: \sim \: 
& \left( \frac{3A_0^2 B_0}{D_0} \right)^{\frac13} \: 
s^{-\frac23} \: t^{\frac13} \:
(\theta^1)^2  \: + \: 
(A_0 B_0 C_0)^{\frac13} \: s^{-1} \: 
(\theta^2)^2 \: + \\
& \left( \frac{B_0 C_0^2 D_0}{3} \right)^{\frac13} \:
s^{-\frac43} \: t^{-\frac13} \: (\theta^3)^2 \: + \: 
\left( \frac{9C_0 D_0}{A_0} \right)^{\frac13} \: 
s^{-\frac13} \: t^{\frac23} \: (\theta^4)^2. \notag
\end{align}

We take coordinates $(x,y,z,r)$ for $G$ in which
$\theta^1 \: = \: dx$, $\theta^2 \: = \: dy \: - \:x dr$,
$\theta^3 \: = \: dz \: - \:  y dr$ and $\theta^4 \: = \: dr$.
Define $\phi_s \: : \: \R^4 \rightarrow \R^4$ by
\begin{equation}
\phi_s(x, y, z, r) \: = \: 
\left( \left( \frac{D_0 s^2}{A_0^2 B_0} \right)^{\frac16}\: x,
\left( \frac{s^3}{A_0 B_0 C_0} \right)^{\frac16} \: y,
\left( \frac{s^4}{B_0 C_0^2 D_0} \right)^{\frac16} \: z,
\left( \frac{A_0 s}{C_0 D_0} \right)^{\frac16} \: r \right).
\end{equation}
Then $\lim_{s \rightarrow \infty} s^{-1} \phi_s^* g(st) \: = \: 
g_\infty(t)$, where
\begin{equation}
g_\infty(t) \: = \: 
3^{\frac13} \: t^{\frac13} \:
(\theta^1)^2  \: + \: 
(\theta^2)^2 \: + \\
3^{-\frac13} \: t^{-\frac13} \: (\theta^3)^2 \: + \: 
3^{\frac23} \: t^{\frac23} \: (\theta^4)^2
\end{equation}
This is an expanding soliton with
$\eta_t(x,y,z,r) \: = \: \left( t^{-\frac26} x, t^{-\frac36} y,
t^{-\frac46} z, t^{-\frac16}r \right)$.

\subsubsection{A7(i)} \label{A7}

The group $G$ is the group $\Sol^4_1$ mentioned above.
The nonzero Lie algebra relations are
$[X_2, X_3] = X_4$, $[X_3, X_1] = X_2$ and $[X_1, X_2] = - X_3$.

The metric is
\begin{equation}
g(t) \: = \: A(t) (\theta^1)^2 \: + \: B(t) (\theta^2)^2 \: + \: 
C(t) (\theta^3)^2 \: + \: D(t) (\theta^4)^2,
\end{equation} 
where
\begin{equation}
d\theta^1 \: = \: 0, \: \: \: \:
d\theta^2 \: = \: - \: \theta^3 \wedge \theta^1, \: \: \: \:
d\theta^3 \: = \: \theta^1 \wedge \theta^2, \: \: \: \:
d\theta^4 \: = \: - \: \theta^2 \wedge \theta^3.
\end{equation}
The Ricci flow asymptotics are
\begin{align}
A(t) \: & \sim \: 4t, \\
B(t) \: & \sim \: (9B_0 C_0 D_0^2)^{\frac16} \: t^{\frac13}, \notag \\ 
C(t) \: & \sim \: (9B_0 C_0 D_0^2)^{\frac16} \: t^{\frac13}, \notag \\ 
D(t) \: & = \: D_0 \: \left(1 + \frac{3D_0}{B_0 C_0} t \right)^{- \frac13}. 
\notag
\end{align}
Then 
\begin{equation}
g_s(t) \: \sim \: 
4t \:
(\theta^1)^2  \: + \: 
(9B_0 C_0 D_0^2)^{\frac16} \: s^{- \frac23} \: t^{\frac13} \:
\left( 
(\theta^2)^2 \: + \: 
(\theta^3)^2 \right) \: + \: \left(
\frac{B_0 C_0 D_0^2}{3} \right)^{\frac13} \: s^{- \frac43} \: 
t^{- \frac13} \: (\theta^4)^2.
\end{equation}

We take coordinates $(x,y,z,r)$ for $G$ in which
$\theta^1 \: = \: dr$, $\theta^2 \: + \: \theta^3 \: = \: dx \: - \:x dr$,
$\theta^2 \: - \: \theta^3 \: = \: dy \: + \:y dr$ and
$\theta^4 \: = \: dz \: +\: \frac14 (xdy - ydx) \: + \: \frac12 \: xydr$.
Define $\phi_s \: : \: \R^4 \rightarrow \R^4$ by
\begin{equation}
\phi_s(x, y, z, r) \: = \: 
\left( (B_0 C_0 D_0^2)^{-\frac{1}{12}} \: s^{\frac13} \: x,
(B_0 C_0 D_0^2)^{-\frac{1}{12}} \: s^{\frac13} \: y,
(B_0 C_0 D_0^2)^{-\frac16} \: s^{\frac23} \: z,
r \right).
\end{equation}
Then $\lim_{s \rightarrow \infty} s^{-1} \phi_s^* g(st) \: = \: 
g_\infty(t)$, where
\begin{equation}
g_\infty(t) \: = \: 
4t \:
(\theta^1)^2  \: + \: 
3^{\frac13} \: t^{\frac13} \:
\left( 
(\theta^2)^2 \: + \: 
(\theta^3)^2 \right) \: + \: 
3^{- \frac13} \: t^{- \frac13} \: (\theta^4)^2.
\end{equation}
This is an expanding soliton with
$\eta_t(x,y,z,r) \: = \: \left( t^{-\frac13} x, t^{-\frac13} y,
t^{-\frac23} z, r \right)$.

\subsubsection{A8} \label{A8}

The group $G$ is a semidirect product $\Nil^3 \widetilde{\times} \R$,
where $\R$ acts on $\Nil^3$, in appropriate coordinates,
by $r \cdot (x,y,z) \: = \:
\left( \cos(r) x + \sin(r) y, - \sin(r) x + \cos(r) y, z \right)$.
The nonzero Lie algebra relations are
$[X_2, X_3] = - X_4$, $[X_3, X_1] = X_2$, $[X_1, X_2] = X_3$.
The $\R$-factor is spanned by $X_1$.

The metric is
\begin{equation}
g(t) \: = \: A(t) (\theta^1)^2 \: + \: B(t) (\theta^2)^2 \: + \: 
C(t) (\theta^3)^2 \: + \: D(t) (\theta^4)^2,
\end{equation} 
where
\begin{equation}
d\theta^1 \: = \: 0, \: \: \: \:
d\theta^2 \: = \: - \: \theta^3 \wedge \theta^1, \: \: \: \:
d\theta^3 \: = \: - \: \theta^1 \wedge \theta^2, \: \: \: \:
d\theta^4 \: = \: \theta^2 \wedge \theta^3.
\end{equation}
The Ricci flow asymptotics are
\begin{align}
A(t) \: & \sim \: \frac{D_0}{2}, \\
B(t) \: & \sim \: (9B_0 C_0 D_0^2)^{\frac16} \: t^{\frac13}, \notag \\ 
C(t) \: & \sim \: (9B_0 C_0 D_0^2)^{\frac16} \: t^{\frac13}, \notag \\ 
D(t) \: & = \: D_0 \: \left(1 + \frac{3D_0}{B_0 C_0} t \right)^{- \frac13}. 
\notag
\end{align}
Then 
\begin{equation}
g_s(t) \: \sim \: 
\frac{D_0}{2} \: s^{-1} \:
(\theta^1)^2  \: + \: 
(9B_0 C_0 D_0^2)^{\frac16} \: s^{- \frac23} \: t^{\frac13} \:
\left( 
(\theta^2)^2 \: + \: 
(\theta^3)^2 \right) \: + \: \left(
\frac{B_0 C_0 D_0^2}{3} \right)^{\frac13} \: s^{- \frac43} \: 
t^{- \frac13} \: (\theta^4)^2.
\end{equation}

Define a diffeomorphism $\phi_s \: : \: \R^4 \rightarrow G$ by
\begin{equation}
\phi_s(x,y,z,r) \: =\: \alpha_s(x,y,z) \: \beta(r), 
\end{equation}
where
$\alpha_s(x,y,z) \: = \: e^{(B_0 C_0 D_0^2)^{- \frac{1}{12}}
s^{\frac13} 
(xX_2+yX_3) + (B_0 C_0 D_0^2)^{-\frac16} s^{\frac23} z X_4}$ and 
$\beta(r) \: = \: e^{D_0^{-\frac12} s^{\frac12} rX_1}$.
Letting $h^{-1} dh$ denote the
Maurer-Cartan form on $G$, we have
\begin{equation}
 \phi_s^* (h^{-1} dh) \: = \: \beta_s^{-1} \: \alpha_s^{-1} \: 
d\alpha_s \: \beta_s \: + \: \beta_s^{-1} \: d\beta_s.
\end{equation}
As conjugation by $\beta_s$ acts isometrically on 
$\alpha_s^{-1} \: d\alpha_s$,
we see that there is a limit as $s \rightarrow \infty$ of
$\phi_s^* g_{s}\left( \cdot \right)$, given by the expanding
soliton $g_\infty(\cdot)$ on $\R \times \Nil^3$; see Section
\ref{nil3}.

\subsubsection{A9} \label{A9} This is a product case
$G \: = \: \widetilde{\SL(2, \R)} \times \R$.

\subsubsection{A10} \label{A10} This is a product case
$G \: = \: \SU(2) \times \R$.

\section{Expanding solitons on vector bundles} \label{vectorbundles}

In this section we consider the expanding soliton equation in the
case of a family $g(\cdot)$ of $\R^N$-invariant metrics
on the total space $M$ of a flat $\R^N$-vector bundle over a manifold
$B$, with the property
that the fiberwise volume forms are preserved by the flat connection.
The vector field $V$ is assumed to be the standard radial vector field
along the fibers.
We show that the expanding soliton equation on $M$ becomes two equations
on B : a harmonic map
equation $G \: : \:  B \rightarrow \SL(N, \R)/\SO(N)$ and an equation
that relates $dG$ to $\Ric_B$. We give examples of expanding soliton
solutions with
$\dim(B) \: = \: 1$, which are generalized $\Sol$-solutions, and an example
with $\dim(B) \: = \: 2$. We then show that if
a rescaling limit $\lim_{k \rightarrow \infty} g_{s_k}(\cdot)$ exists for
such an $\R^N$-invariant Ricci flow then the limit satisfies the
harmonic-Einstein equations.

The expanding solitons $(M_\infty, g_\infty(t))$ of 
Section \ref{homogeneous} all have a certain
fibration structure.  Namely, there is a Riemannian submersion
$\pi \: : \: M_\infty \rightarrow B_\infty$ whose fibers
are diffeomorphic to a nilpotent Lie group ${\mathcal N}$
and whose holonomy preserves the natural flat linear connection
$\nabla^{aff}$
on ${\mathcal N}$. (The connection $\nabla^{aff}$ has the property
that left-invariant vector fields are parallel.)
The diffeomorphisms $\{\eta_t\}_{t > 0}$ act fiberwise and
arise from a $1$-parameter group $\{a_t\}_{t>0}$ of automorphisms
of ${\mathcal N}$, by $\eta_t \: = \: a_{t^{-1}}$. 
We list below the relevant groups ${\mathcal N}$ that appeared in
Section \ref{homogeneous}, along with the subsection in which they 
appeared.
(There are also some product cases that we omit.)
\begin{align}
\{e\} & \: \: \: \: \: \: \ref{2d},\ref{H3},\ref{h4},\ref{ch2}\\
\R & \: \: \: \: \: \: \ref{1d},\ref{sl2} \notag \\
\R^2 & \: \: \: \: \: \: \ref{sol},\ref{f4},\ref{A9} \notag \\
\R^3 & \: \: \: \: \: \: \ref{isomr2},\ref{A2},\ref{A3},\ref{A5} \notag \\
\R^4 & \: \: \: \: \: \: \ref{A1} \notag \\
\Nil^3 & \: \: \: \: \: \: \ref{nil3},\ref{A4},\ref{A7},\ref{A8} \notag \\
\Nil^4 & \: \: \: \: \: \: \ref{A6} \notag
\end{align}

These special fibration structures are related to
the results of Cheeger-Fukaya-Gromov on collapsing with
bounded sectional curvature 
\cite{Cheeger-Fukaya-Gromov (1992)}. Namely, any sufficiently 
collapsed manifold can be slightly perturbed to have a so-called
$\Nil$-structure, where ``$\Nil$'' refers to a local nilpotent Lie
algebra of Killing vector fields. It will follow from Section \ref{discussion} 
that if
$(M, p, g(\cdot))$ is a pointed type-III Ricci flow solution then there is
a sequence $\{s_j\}_{j=1}^\infty$ tending to infinity so that
there is a limit flow
$g_\infty(\cdot) \: = \: \lim_{j \rightarrow \infty} g_{s_j}(\cdot)$ in an 
appropriate sense which, if
$\lim_{t \rightarrow \infty} t^{- \frac12} \inj_{g(t)}(p) \: = \: 0$, will have a $\Nil$-structure. 

This suggests looking for expanding soliton solutions with a special
fibration structure of the type mentioned above,
with the action of the 
diffeomorphisms $\{\eta_t\}_{t > 0}$ being compatible with
the fibration structure. In terms of the fiber, it is known that
there are many nilpotent Lie groups that do admit Ricci soliton
metrics and also some that do not 
\cite{Lauret (2001)}.

Based on these considerations, in this section we look at what the
expanding soliton equation becomes if we assume compatibility with
the simplest type of $\Nil$-structure.
Namely, we consider the expanding soliton equation on the total
space of an $\R^N$-vector bundle $\pi \: : \: M \rightarrow B$. We assume
that \\
1. $\pi$ is a Riemannian submersion. \\
2. For each $b \in B$, there is a neighborhood $U_b$ of $b$ in $B$ so that
there is a
free isometric $\R^N$-action on $\pi^{-1}(U_b)$ which acts by 
translation on the fibers.\\
3. The diffeomorphism $\eta_t$ is
fiberwise multiplication by $t^{-\frac12}$. \\

Let $s \: : \:  B \rightarrow M$ be the zero-section.
The $\R^N$-action on $\pi^{-1}(U_b)$ gives a local trivialization
of the vector bundle by
$(b^\prime, \vec{v}) \rightarrow s(b^\prime) + \vec{v}$.
In view of this, it is natural to reduce the data to \\
1. A vector bundle on $B$ with a flat vector bundle connection $\nabla$,\\
2. A Riemannian metric on $B$ and\\
3. Flat Riemannian metrics on the fibers. \\

There is a corresponding canonical metric on $M$.
We do not assume that $\nabla$ preserves the fiberwise metrics.

For simplicity of notation, in this section we write 
$\overline{g}_{IJ}$ for the
metric on $M$, $\overline{R}_{IJ}$ for the Ricci tensor on $M$,
etc. We write $g_{\alpha \beta}$ for the
metric on $B$, $R_{\alpha \beta}$ for the Ricci tensor on $B$,
etc. We let Greek indices
denote horizontal directions and we let lower case Roman indices denote
vertical directions. 
In terms of local coordinates 
$\{x^\alpha, x^i\}$, we can write the metric on 
$M$ as
\begin{align}
\overline{g}_{\alpha \beta} \: & = \: g_{\alpha \beta}(b) \\
\overline{g}_{i\alpha} \: & = \: 0 \notag \\
\overline{g}_{ij} \: & = \: g_{ij}(b) \notag.
\end{align}
We will use the Einstein summation convention freely.

Hereafter
we assume that $\nabla$ preserves the fiberwise volume forms, as
this is what arises in the examples of Section \ref{homogeneous}.
We write
\begin{equation}
g_{ij;\alpha \beta} \: = \: g_{ij,\alpha \beta} \: - \:
\Gamma^\sigma_{\: \: \alpha \beta} \: g_{ij, \sigma}.
\end{equation}

\begin{proposition} \label{harmoniceinstein}
The expanding soliton equation becomes the pair of equations
\begin{equation} \label{baseflow}
R_{\alpha \beta} \: - \: \frac14 \:
g^{ij} \: g_{jk, \alpha} \: g^{kl} \: g_{li, \beta} \: + \:
\frac{1}{2t} \: g_{\alpha \beta} \: = \: 0
\end{equation}
and
\begin{equation} \label{harmonic}
g^{\alpha \beta} \:
g_{ij; \alpha \beta} \: - \: g^{\alpha \beta} \:
g_{ik, \alpha} \: g^{kl} \: g_{lj, \beta} \: = \: 0.
\end{equation}
\end{proposition}
\begin{proof}
The nonzero Christoffel symbols are
\begin{align}
\overline{\Gamma}^\alpha_{\: \: \beta \gamma} \: & = \: 
{\Gamma}^\alpha_{\: \: \beta \gamma} \\
\overline{\Gamma}^\alpha_{\: \: ij} \: & = \: - \: \frac12 \:
g^{\alpha \beta} \: g_{ij, \beta} \notag \\ 
\overline{\Gamma}^i_{\: \: j \alpha} \: = \:
\overline{\Gamma}^i_{\: \: \alpha j} \: & = \: \frac12 \: 
g^{ik} \: g_{kj, \alpha}. \notag 
\end{align}
The nonzero components of the curvature tensor are
\begin{align} \label{curvature}
\overline{R}^\alpha_{\: \: \beta \gamma \delta} \: & = \:
{R}^\alpha_{\: \: \beta \gamma \delta} \\
\overline{R}^\alpha_{\: \: i \beta j} \: & = \: - \: \frac12 \:
g^{\alpha \gamma} \: g_{ij;\gamma \beta} \: + \: \frac14 \:
g^{\alpha \gamma} \: g_{ik,\beta} \: g^{kl} \: g_{lj, \gamma} \notag \\
\overline{R}^i_{\: \: j \alpha \beta} \: & = \: - \: \frac14 \:
g^{ik} \: g_{kl, \alpha} \: g^{lm} \: g_{mj, \beta} \: + \: \frac14 \:
g^{ik} \: g_{kl, \beta} \: g^{lm} \: g_{mj, \alpha} \notag \\
\overline{R}^i_{\: \: jkl} \: & 
= \: - \: \frac14 \: g^{ir} \: g_{rk, \alpha} \:
g^{\alpha \beta} \: g_{jl, \beta} \: + \: 
\frac14 \: g^{ir} \: g_{rl, \alpha} \:
g^{\alpha \beta} \: g_{jk, \beta}. \notag
\end{align}
The Ricci tensor is
\begin{align} \label{ricci}
\overline{R}_{\alpha \beta} \: & = \: R_{\alpha \beta} \: - \:
\frac12 \: g^{ij} \: g_{ij;\alpha \beta} \: + \: \frac14 \:
g^{ij} \: g_{jk, \alpha} \: g^{kl} \: g_{li, \beta} \\
\overline{R}_{\alpha i} \: & = \: 0 \notag \\
\overline{R}_{ij} \: & = \: - \: \frac12 \: g^{\alpha \beta} \:
g_{ij; \alpha \beta} \: + \: \frac12 \: g^{\alpha \beta} \:
g_{ik, \alpha} \: g^{kl} \: g_{lj, \beta}
\: - \: \frac14 \: g^{\alpha \beta} \: g^{kl} \: g_{kl, \alpha}
\: g_{ij, \beta}. \notag
\end{align}

As $\nabla$ preserves the fiberwise volume forms,
\begin{equation}
g^{ij} \: g_{ij,\alpha} \: = \: 0
\end{equation}
and
\begin{equation}
g^{ij} \: g_{ij; \alpha \beta} \:  = \: g^{ij} \: g_{jk,\beta} \:
g^{kl} \: g_{li, \alpha}.
\end{equation}
Then
\begin{align} \label{ricci2}
\overline{R}_{\alpha \beta} \: & = \: R_{\alpha \beta} \: - \: \frac14 \:
g^{ij} \: g_{jk, \alpha} \: g^{kl} \: g_{li, \beta} \\
\overline{R}_{\alpha i} \: & = \: 0 \notag \\
\overline{R}_{ij} \: & = \: - \: \frac12 \: g^{\alpha \beta} \:
g_{ij; \alpha \beta} \: + \: \frac12 \: g^{\alpha \beta} \:
g_{ik, \alpha} \: g^{kl} \: g_{lj, \beta}. \notag
\end{align}

If $\{V(t)\}_{t > 0}$ is the vector field 
generating $\{\eta_t\}_{t > 0}$ then
\begin{align}
({\mathcal L}_V g)_{\alpha \beta} \: & = \: 0 \\
({\mathcal L}_V g)_{\alpha i} \: & = \: 0 \notag \\
({\mathcal L}_V g)_{ij} \: & = \: - \: \frac{1}{t} \: g_{ij}. \notag
\end{align}
The proposition follows. 
\end{proof}

\begin{remark}
Without assuming that $\nabla$ preserves the fiberwise volume forms,
it follows from (\ref{ricci}) that
$g^{ij} \: \overline{R}_{ij} \: = \: - \:
\frac{1}{\sqrt{|G|}} \: g^{\alpha \beta} \: 
\left( \sqrt{|G|} \right)_{; \alpha \beta}$,
where $|G| \: = \: \det \left( g_{ij} \right)$.
Then under the Ricci flow,
\begin{equation}
\frac{\partial \sqrt{|G|}}{\partial t} \: = \:
- \: \sqrt{|G|} \: g^{ij} \: \overline{R}_{ij} \: = \:
g^{\alpha \beta} \: \left( \sqrt{|G|} \right)_{; \alpha \beta}.
\end{equation}
Hence it is consistent to assume that $|G|$ is spatially constant,
i.e. that $\nabla$ preserves the fiberwise volume forms. \\
\end{remark}

We will call equations (\ref{baseflow})-(\ref{harmonic}) the 
{\em harmonic-Einstein} equations.

Now consider the space ${\mathcal S}$ of 
positive-definite symmetric  matrices 
$\{G_{ij}\}$ on $\R^N$ with a fixed determinant.
An element $A \in \SL(N, \R)$ acts on ${\mathcal S}$ by sending
$G$ to $A G A^T$. This identifies ${\mathcal S}$ with
$\SL(N, \R)/\SO(N)$.
The corresponding Riemannian metric on ${\mathcal S}$ can be
written informally as $\Tr \left( G^{-1} dG \right)^2$. That is,
for a symmetric matrix $K \in T_G {\mathcal S}$
\begin{equation}
\langle K, K \rangle_G \: = \: G^{ij} \: K_{jk} \: G^{kl} \: K_{li}.
\end{equation}

\begin{proposition}
Equation (\ref{harmonic}) is the local expression for a harmonic map from
$B$ to ${\mathcal S}$.
\end{proposition}
\begin{proof}
The energy of a map $G \: : \: B \rightarrow {\mathcal S}$ is
\begin{equation}
E(G) \: = \: \frac12 \:
\int_B g^{\alpha \beta} \: \Tr \left(  G^{-1} \: G,_\alpha \:
G^{-1} \: G,_\beta \right) \: \dvol. 
\end{equation}
Consider a variation of $G$ of the form $\delta G \: = \: KG$
with $\Tr K = \: 0$. The variation of $E$ is
\begin{equation}
\delta E \: = \: \int_B g^{\alpha \beta} \: 
\Tr \left(  G^{-1} \: G_{,\alpha} \:
G^{-1} \: K_{,\beta} \: G \right) \: \dvol
 \: = \: \int_B g^{\alpha \beta} \: 
\Tr \left(G_{,\alpha} \:
G^{-1} \: K_{,\beta} \right) \: \dvol.
\end{equation}
If $K$ is compactly supported then integration by parts gives
\begin{equation}
\delta E \: = \: - \: \int_B \: g^{\alpha \beta} \: 
\Tr \left(\left( G_{,\alpha} \: G^{-1} \right)_{;\beta} \: K \right)
\: \dvol.
\end{equation}
If this vanishes for all such $K$ then
\begin{equation} \label{tracing}
g^{\alpha \beta} \:
\left( G_{,\alpha} \: G^{-1} \right)_{;\beta} \:  = \: \sigma \: I
\end{equation}
for some function $\sigma$ on $B$.
On the other hand, as $G$ has constant determinant, 
\begin{equation}
\Tr \left( G_{,\alpha} \: G^{-1} \right) \: = \: 0
\end{equation}
and so
\begin{equation}
\Tr \left( G_{,\alpha} \: G^{-1} \right)_{;\beta} \: = \: 0.
\end{equation}
Tracing (\ref{tracing}) gives $\sigma \: = \: 0$, so the
variational equation is
\begin{equation}
g^{\alpha \beta} \: \left( G_{,\alpha} \: G^{-1} \right)_{;\beta} \:  = \: 
g^{\alpha \beta} \: G_{;\alpha \beta} \: G^{-1}
- \: g^{\alpha \beta} \: G_{,\alpha} \: G^{-1} \: 
G_{,\beta} \: G^{-1} \: = \: 0.
\end{equation}
Equivalently,
\begin{equation}
g^{\alpha \beta} \: G_{;\alpha \beta} \: 
- \: g^{\alpha \beta} \: G_{,\alpha} \: G^{-1} \: 
G_{,\beta} \: = \: 0,
\end{equation}
which is the same as (\ref{harmonic}).
\end{proof}

The equations (\ref{baseflow}) and (\ref{harmonic}) are defined locally.
To give them global meaning, 
let $\rho \: : \: 
\pi_1(B) \rightarrow \SL(N, \R)$ be the holonomy representation of the
flat connection $\nabla$.
Let $\widetilde{B}$ be the universal cover
of $B$. The flat vector bundle $M$ over $B$ can be written as
$\widetilde{B} \times_{\pi_1(B)} \R^N$, where 
$\pi_1(B)$ is represented on $\R^N$ via $\rho$.
Then the family of fiberwise metrics $\{g_{ij}(b) \}_{b \in B}$ corresponds 
to a $\pi_1(B)$-equivariant map $G \: : \: \widetilde{B} \rightarrow 
\SL(N, \R)/\SO(N)$, where $\pi_1(B)$ acts on $\SL(N, \R)/\SO(N)$
via left multiplication by $\rho$. Equation (\ref{harmonic})  says that
$G$ is a harmonic map.  As $\overline{R}_{ij} \: = \: 0$, the metrics 
$\{g_{ij}(b) \}_{b \in B}$ are constant in time and so
the map $G$ is time-independent. The metric $g_{\alpha \beta}$ on 
${B}$ is proportionate to $t$. Equation (\ref{baseflow}) relates
the metric on $B$ to the harmonic map $G$.

If, after an appropriate choice of basis, the representation $\rho$ takes
value in $\SL(N, \Z)$ then we can quotient $M$ fiberwise by
$\Z^N$. The resulting space 
$\widetilde{M}/(\Z^N \widetilde{\times} \pi_1(B))$ is a flat torus bundle
over $B$ and carries a quotient Ricci flow metric.

\begin{example}
If $B$ is compact and $\rho$ is trivial then $G$ must be a point map
(see, for example, \cite[Section 1.2]{Jost-Zuo (1997)}). Thus 
the expanding soliton metric $(M, \overline{g}(t))$ is just a product metric
$(B, t g_B) \times (\R^N, g_{flat})$ for an Einstein metric $g_B$ on $B$
with Einstein constant $- \: \frac12$.
\end{example}

\begin{example}
Suppose that $B \: = \: \R$. Let $X$ be a real diagonal
$(N \times N)$-matrix with vanishing trace.  Then there is a geodesic
$G \: : \: \R \rightarrow {\mathcal S}$ given by
$G(s) \: = \: e^{sX}$. Equation (\ref{baseflow}) is satisfied 
by the metric $\frac{t}{2} \: \Tr(X^2) \: ds^2$ on $B$. This
generalizes the $\Sol$-solutions of Sections \ref{sol} and
\ref{A2}.

If $e^{\frac{X}{2}}$ is conjugate to an integer matrix $A \in \SL(N, \Z)$ then
one obtains a quotient Ricci flow solution on the total space of a
flat $T^N$-bundle over $S^1$ with holonomy $A$.
\end{example}

\begin{example} \label{F4}
Take $B \: = \: \SL(N, \R)/\SO(N)$. The identity map $G$ from $B$ to
$\SL(N, \R)/\SO(N)$ is harmonic. 
Equation (\ref{baseflow}) is satisfied by the
canonical metric on $B$, after appropriate normalization.  
As $B$ is the moduli space for constant-volume inner products on 
$\R^N$, the manifold $M$ is the corresponding universal Euclidean
bundle over $B$. 
If $\Gamma$ is a finite-index torsion-free subgroup of $\SL(N, \Z)$ then
there is a quotient $T^n$-bundle over $\Gamma \backslash \SL(N, \Z)/\SO(N)$,
with a quotient Ricci flow solution.

In the case $N = 2$, we can identify $\SL(2, \R)/\SO(2)$ with
the space of translation-invariant complex structures on $\R^2$.
The manifold  $M$ is the resulting universal complex line bundle. It has a
homogeneous complex geometry, called F4 in \cite{Wall (1986)}.
The quotient of $M$ by $\Z^2 \widetilde{\times} \SL(2, \Z)$ is the universal curve of
complex structures on $T^2$, which we equip with constant volume. 
\end{example}

Turning from the expanding soliton equation, 
we now consider Ricci flow on a $1$-parameter family of
metrics $\overline{g}(\cdot)$ on $M$ of the
type considered above. We ask for sufficient conditions to ensure
that the limit of a convergent subsequence
of rescaled flows $s_j^{-1} \: \overline{g}(s_j \cdot)$
(modulo diffeomorphisms) is in fact
an expanding soliton. This is a different question than that addressed in
Proposition \ref{limit} where we assumed that there is an actual limit
as $s \rightarrow \infty$ of $s^{-1} \: g(s \cdot)$, instead of just
a convergent subsequence.

Suppose first that $N = 0$, i.e. we just have Ricci flow 
on $B$.
Given a type-III Ricci flow $g(\cdot)$ on a compact manifold $B$, 
suppose that there is a sequence $\{s_j\}_{j=1}^\infty$ of positive numbers
converging to infinity and diffeomorphisms
$\{\phi_j\}_{j=1}^\infty$ of $B$ so that
$\left\{ s_j^{-1} \phi_j^* g(s_j \cdot) \right\}_{j=1}^\infty$ converges
to a Ricci flow $g_\infty(\cdot)$ on $B$.
Then $g_\infty(t) \: = \: t \: g_E$ for a Einstein metric
$g_E$ on $B$ satisfying $\Ric(g_E) \: = \: - \: \frac12 \: g_E$
\cite[Theorem 1.3]{Feldman-Ilmanen-Ni (2005)},\cite[Section 7]{Hamilton 
(1999)}.
The hypotheses are equivalent to saying that we 
have a type-III Ricci flow solution $g(\cdot)$ 
on a manifold $B$ such that  $\limsup_{t \rightarrow \infty} \: t^{-\frac12} \:
\diam(B, g(t)) \: < \: \infty$ and $\liminf_{t \rightarrow \infty} \: 
t^{-\frac12} \:
\inj(B, g(t)) \: > \: 0$.

We now consider the case $N > 0$. An automorphism of a flat
vector bundle $(W, \nabla^{flat})$ 
over $B$ is an invertible vector bundle map 
$\widehat{\phi} \: : \: W \rightarrow W$ with $\widehat{\phi}^* 
\nabla^{flat}
\: = \: \nabla^{flat}$. It covers a diffeomorphism $\phi$ of $B$.

\begin{proposition} \label{helimits}
Fix a flat $\R^N$-vector bundle $M$ on a compact manifold $B$.
Let $\{\overline{g}(t)\}_{t \in (0, \infty)}$ 
be a $1$-parameter family of metrics on $M$
of the type considered above. Suppose that 
$\overline{g}(\cdot)$ satisfies the Ricci flow equation.
Suppose that\\
1. There is a sequence $\{s_j\}_{j=1}^\infty$ of positive numbers
converging to infinity and\\
2. There are automorphisms
$\{\widehat{\phi}_j\}_{j=1}^\infty$ of the flat vector bundle $M$ so that \\
3. $\left\{ s_j^{-1} \widehat{\phi}_j^* \overline{g}(s_j \cdot)
\right\}_{j=1}^\infty$ converges
to a Ricci flow solution $\overline{g}_\infty(\cdot)$.

Then $\overline{g}_\infty(\cdot)$ satisfies the harmonic-Einstein equations
(\ref{baseflow})-(\ref{harmonic}).
\end{proposition}
\begin{proof}
The proof is an adaptation of one of the proofs of
\cite[Theorem 1.3]{Feldman-Ilmanen-Ni (2005)}, this particular
proof being due to Hamilton.
We write $G_{ij}$ for $g_{ij}$. On $M$ we have the equation
\begin{equation}
\frac{\partial \overline{R}}{\partial t} \: = \:
\overline{\triangle} \: \overline{R} \: + \: 2 \:
|\overline{R}_{IJ}|^2.
\end{equation}
Translating this to an equation on $B$ gives
\begin{align} \label{maximum}
& \frac{\partial}{\partial t} \left( R \: - \:
\frac14 \: g^{\alpha \beta} \: \Tr \left( G^{-1} \: G_{, \alpha} \:
G^{-1} \: G_{, \beta} \right) \right) \: = \: 
\triangle \left( R \: - \:
\frac14 \: g^{\alpha \beta} \: \Tr \left( G^{-1} \: G_{, \alpha} \:
G^{-1} \: G_{, \beta} \right) \right)  + \\
& 2 \: \left| R_{\alpha \beta} \: - \:
\frac14 \: \Tr \left( G^{-1} \: G_{, \alpha} \:
G^{-1} \: G_{, \beta} \right) \right|^2 \: + \notag \\
& 2 \: \Tr \left( - \: \frac12 \: g^{\alpha \beta} \: G^{-1} \: 
G_{;\alpha \beta} \: + \: \frac12 \: g^{\alpha \beta} \:
G^{-1} \: G_{, \alpha} \: G^{-1} \: G_{, \beta} \right)^2.\notag 
\end{align}
The maximum principle implies that
\begin{equation}
R \: - \:
\frac14 \: g^{\alpha \beta} \: \Tr \left( G^{-1} \: G_{, \alpha} \:
G^{-1} \: G_{, \beta} \right) \: + \: \frac{n}{2t} \: \ge \: 0,
\end{equation}
where $n = \dim(B)$.
On the other hand, putting $\tilde{V}(t) \: = \: t^{- \frac{n}{2}} \: 
\int_B \dvol_B$, one has
\begin{equation} \label{mono}
\frac{d\tilde{V}}{dt} \: = \: - \: t^{- \frac{n}{2}} \: \int_B
\left( R \: - \:
\frac14 \: g^{\alpha \beta} \: \Tr \left( G^{-1} \: G_{, \alpha} \:
G^{-1} \: G_{, \beta} \right) \: + \: \frac{n}{2t} \right) \: \dvol_B.
\end{equation}
Then $\tilde{V}(t)$ is nonincreasing in $t$ and
the corresponding quantity $\tilde{V}_\infty(t)$ for
$\overline{g}_\infty(\cdot)$ must be constant in $t$.
Equation (\ref{mono}), applied to $\overline{g}_\infty(\cdot)$, implies
\begin{equation}
R_\infty \: - \:
\frac14 \: g_\infty^{\alpha \beta} \: \Tr \left( G_\infty^{-1} \: 
G_{\infty, \alpha} \:
G_\infty^{-1} \: G_{\infty, \beta} \right) \: + \: \frac{n}{2t} \: = \: 0.
\end{equation}
Plugging this into (\ref{maximum}) (applied to $\overline{g}_\infty$)
gives
\begin{align}
\frac{n}{2t^2} \: = \: & 
2 \: \left| R_{(\infty)\alpha \beta} \: - \:
\frac14 \: \Tr \left( G_\infty^{-1} \: G_{\infty, \alpha} \:
G_\infty^{-1} \: G_{\infty, \beta} \right) \right|^2 \: + \\
& 2 \: \Tr \left( - \: \frac12 \: g_\infty^{\alpha \beta} \: 
G_\infty^{-1} \: 
G_{\infty;\alpha \beta} \: + \: \frac12 \: g_\infty^{\alpha \beta} \:
G_\infty^{-1} \: G_{\infty, \alpha} \: G_\infty^{-1} \: 
G_{_\infty, \beta} \right)^2,\notag 
\end{align}
or
\begin{align}
0 \: = \: & 
2 \: \left| R_{(\infty)\alpha \beta} \: - \:
\frac14 \: \Tr \left( G_\infty^{-1} \: G_{\infty, \alpha} \:
G_\infty^{-1} \: G_{\infty, \beta} \right) 
\: + \: \frac{1}{2t} \: g_{\infty, \alpha \beta} \right|^2 \: + \\
& 2 \: \Tr \left( - \: \frac12 \: g_\infty^{\alpha \beta} \: 
G_\infty^{-1} \: 
G_{\infty;\alpha \beta} \: + \: \frac12 \: g_\infty^{\alpha \beta} \:
G_\infty^{-1} \: G_{\infty, \alpha} \: G_\infty^{-1} \: 
G_{_\infty, \beta} \right)^2. \notag 
\end{align}
The proposition follows.
\end{proof}

Based on Proposition \ref{helimits}, one can speculate the expanding solitons
that arise in type-III Ricci flows again involve some sort of harmonic-Einstein equations.
Of course, the relevant nilpotent Lie groups may be more complicated than 
$\R^N$ and they may act with orbits of varying dimensions. Nevertheless, a 
$\Nil$-structure has a quotient space with bounded geometry in a certain sense
\cite[Appendix 1]{Cheeger-Fukaya-Gromov (1992)}, which is where the
harmonic-Einstein equations would live.

\section{Ricci Flow on \'Etale Groupoids} \label{rgroupoids}

An \'etale groupoid is a mathematical object which in some sense
combines the notions of topological spaces and discrete groups.
A Riemannian groupoid is an \'etale groupoid equipped with an
invariant Riemannian metric.
The relevance
for us comes from the Cheeger-Fukaya-Gromov theory of bounded
curvature collapse, which implies that when a Riemannian manifold
collapses with bounded sectional curvature, 
it asymptotically obtains local symmetries.

In this section we recall some basic definitions about Riemannian groupoids.
A good source for background information is 
\cite[Section 7]{Haefliger (2001)}. Further references are
\cite[Chapter IIIG]{Bridson-Haefliger (1999)},
\cite{Haefliger (1988)} and \cite[Appendix D]{Salem (1988)}.  We then
prove an extension of Hamilton's compactness theorem,  not
assuming a lower bound on the injectivity radius. Although it takes a bit of
time to set up the right framework, once the framework is in place then
the proof is almost the same as in Hamilton's paper \cite{Hamilton (1995)}.
We discuss how the long-time behavior of type-III Ricci flow solutions
becomes a problem about the dynamics of the $\R^+$-action on a 
compact space ${\mathcal S}_{n,K}$.  We list the
Riemannian groupoids that arise in the long-time
behavior of finite-volume locally homogeneous three-dimensional Ricci 
flow solutions.

\subsection{Etale groupoids}

A {\em groupoid} $G$ consists of \\
1. Sets $G^{(0)}$ and $G^{(1)}$,\\
2. An injection
$e \: : \: G^{(0)} \rightarrow G^{(1)}$ (with which we will think of 
the ``units'' $G^{(0)}$ as a subset
of $G^{(1)}$), \\
3. ``Source'' and ``range'' 
maps $s, r \: : \: G^{(1)} \rightarrow G^{(0)}$ with
$s\circ e \: = \: r\circ e \: = \: \Id \big|_{G^{(0)}}$ and \\
4. A partially-defined
multiplication $G^{(1)} \times G^{(1)} \rightarrow G^{(1)}$

so that\\
1. The product $\gamma \gamma^\prime$ of $\gamma, \gamma^\prime \in G^{(1)}$
is defined if and only if $s(\gamma) \: = \: r(\gamma^\prime)$, and then
$s(\gamma \gamma^\prime) \: = \: s(\gamma^\prime)$ and 
$r(\gamma \gamma^\prime) \: = \: r(\gamma)$. \\
2. $(\gamma \gamma^\prime) \gamma^{\prime \prime} \: = \: \gamma 
(\gamma^\prime \gamma^{\prime \prime})$
whenever the two sides are defined. \\
3. $\gamma s(\gamma) \: = \: r(\gamma) \gamma \: = \: \gamma$. \\
4. For all $\gamma \in G^{(1)}$, there is an element 
$\gamma^{-1} \in G^{(1)}$ so that
$\gamma \gamma^{-1} \: = \: r(\gamma)$ and $\gamma^{-1} \gamma \: = \: 
s(\gamma)$.\\

A {\em morphism} $m \: : \: G_1 \rightarrow G_2$
between two groupoids is given by maps
$m^{(1)} \: : \: G_1^{(1)} \rightarrow G_2^{(1)}$ and 
$m^{(0)} \: : \: G_1^{(0)} \rightarrow G_2^{(0)}$ that satisfy obvious
compatibility conditions. An {\em isomorphism} between $G_1$ and $G_2$
is an invertible morphism. 

Given $x \in G^{(0)}$, 
we write $G^x \: = \: r^{-1}(x)$, $G_x \: = \: s^{-1}(x)$ and
$G^x_x \: = \: G^x \cap G_x$, the latter being the
{\em isotropy group} of $x$. The {\em orbit} of $x$ is the
set $O_x \: =  \: s(r^{-1}(x))$. 
There is an {\em orbit space} ${\mathcal O}$.

A {\em pointed} groupoid $(G, O_x)$ is a groupoid $G$ equipped with a 
preferred
orbit $O_x$. A morphism $m \: : \: (G_1, O_{x_1}) \rightarrow 
(G_2, O_{x_2})$ of pointed groupoids will be assumed to have the property
that $m^{(0)}$ sends $O_{x_1}$ to $O_{x_2}$.

A groupoid $G$ is {\em smooth} if 
$G^{(1)}$ and $G^{(0)}$ are smooth manifolds, 
$e$ is a smooth embedding, $s$ and $r$ are
submersions, and the structure maps are all
smooth. (For example, multiplication is supposed to be a smooth map
from $\{(\gamma, \gamma^\prime) \in G^{(1)} \times G^{(1)} 
\: : \: s(\gamma) \: = \: r(\gamma^\prime)\}$ to $G^{(1)}$.)
Morphisms between smooth groupoids are assumed to be smooth.
A smooth groupoid is
{\em \'etale} if $s$ and $r$ are local diffeomorphisms.
Hereafter, unless otherwise stated, groupoids will be smooth and
\'etale. We do not assume
that $G^{(1)}$ has a countable basis, although in the cases of interest
$G^{(0)}$ will have a countable basis.

An \'etale groupoid $G$ is {\em Hausdorff} if
$G^{(1)}$ is Hausdorff and whenever $c \: : \: [0, 1) \rightarrow
G^{(1)}$ is a continuous path such that $\lim_{t \rightarrow 1} s(c(t))$ and
$\lim_{t \rightarrow 1} r(c(t))$ exist, there is a limit
$\lim_{t \rightarrow 1} c(t)$ in $G^{(1)}$. Hereafter we assume that
$G$ is Hausdorff.  

\begin{example} \label{example1}
If $M$ is a smooth manifold then there is an \'etale groupoid $G$ with
$G^{(1)} \: = \: G^{(0)} \: = \: M$, where
$s$ and $r$ are the identity maps. The product $m \cdot m^\prime$ is
defined if and only if $m = m^\prime$, in which case the product
is $m$. We will call
this the groupoid $M$.
\end{example}

\begin{example} \label{example2}
If a discrete 
group $\Gamma$ acts smoothly on $M$, define the cross-product groupoid
$G \: = \: M \rtimes \Gamma$ as follows. Put 
$G^{(1)} \: = \: M \times \Gamma$ and $G^{(0)} \: = \: M$, with
$r(m, \gamma) \: = \: m$ and $s(m, \gamma) \: = \: m \gamma$.
The product $(m, \gamma) \cdot (m^\prime, \gamma^\prime)$ is defined if
$m\gamma = m^\prime$, in which case the product is $(m, \gamma \gamma^\prime)$.

For example, we can take $\Z$ acting on $S^1$ by rotations.
Or we can take $\SO(2)$ acting on $S^1$ by rotations. Note that in the latter
case we give $\SO(2)$ the discrete topology so that $G$ will be \'etale.
\end{example}

\subsection{Equivalence of \'etale groupoids}

If ${\mathcal U} \: = \: 
\{U_i\}_{i \in I}$ is an open cover of $G^{(0)}$ then there is
a new \'etale groupoid $G_{\mathcal U}$, called the {\em localization} of $G$
to ${\mathcal U}$.
It has $G_{\mathcal U}^{(1)} \: = \: \{(i,\gamma,j) \: : \: 
s(\gamma) \in U_j, r(\gamma)
\in U_i\}$ and $G_{\mathcal U}^{(0)} \: = \: \{(i, \gamma, i) \: : \: 
\gamma \in U_i\}$.
The source and range maps send
$(i,\gamma,j)$ to $(j, s(\gamma), j)$ and $(i, r(\gamma), i)$, 
respectively. The product
is $(i,\gamma,j) \cdot (j, \gamma^\prime, k) \: = \: 
(i, \gamma \gamma^\prime, k)$. 

Two \'etale groupoids $G_1$ and $G_2$ are {\em equivalent} if there are
open covers ${\mathcal U}_1$ and ${\mathcal U}_2$ of 
$G_1^{(0)}$ and $G_2^{(0)}$, respectively, so that
the localizations $(G_1)_{{\mathcal U}_1}$ and $(G_2)_{{\mathcal U}_2}$
are isomorphic. Other ways of expressing this are given in
\cite[p. 597, 599-601]{Bridson-Haefliger (1999)} but
the above definition is good enough for our purposes. 

\begin{example} \label{example3}
Let $G = M$ be the groupoid of Example \ref{example1}. Let
$\{U_i\}_{i=1}^\infty$ be an open cover of $M$. Then $G_{\mathcal U}^{(1)}$
is the disjoint union of the $U_i$'s and $G_{\mathcal U}^{(1)}$
consists of pairs of points $(m_i, m_j) \in U_i \times U_j$ that
get identified to the same point in $M$.  By definition,
$G_{\mathcal U}$ is equivalent to $G = M$.
\end{example}

\begin{example} \label{example4}
In the setup of Example \ref{example2}, suppose that $\Gamma$ acts freely
and properly discontinuously on $M$. Then the cross-product groupoid
$M \rtimes \Gamma$ is equivalent to the groupoid $M/\Gamma$.
\end{example}

\begin{remark}
We wish to identify equivalent groupoids. A more intrinsic approach
is to consider instead the category ${\mathcal B}G$ of 
$G$-sheaves. Here a $G$-sheaf
is a local homeomorphism $\sigma \: : \:
E \rightarrow G^{(0)}$ equipped with a continuous
right $G$-action $E \times_{G^{(0)}} {G^{(1)}} \rightarrow E$.
Pulling back the differentiable structure 
to $E$, we may assume that
$\sigma$ is a local diffeomorphism. Then ${\mathcal B}G$ is a topos and
equivalent groupoids give rise to equivalent topoi
\cite{Maclane-Moerdijk (1996),Moerdijk (1995)}. However, we will not
pursue this approach.
\end{remark}

\subsection{Riemannian groupoids} \label{Riemannian groupoids}

A {\em smooth path} $c$ in $G$ consists of a partition 
$ 0 = t_0 \le t_1 \le \ldots \le t_k = 1$, and a sequence
$c \: = \: (\gamma_0, c_1, \gamma_1, \ldots, c_k, \gamma_k)$ where
$c_i\: : \: [t_{i-1}, t_i] \rightarrow G^{(0)}$ is smooth, 
$\gamma_i \in G^{(1)}$,
$c_{i}(t_{i-1}) \: = \: s(\gamma_{i-1})$ and  
$c_{i}(t_{i}) \: = \: r(\gamma_{i})$. 
The path is said to go from $r(\gamma_0)$ to
$s(\gamma_k)$.
The groupoid $G$ is {\em path connected} if any $x, y \in G^{(0)}$ can be
joined by a smooth path. We will generally assume that $G$ is path
connected.

Given $\gamma \in G^{(1)}$, there is a neighborhood $U$ of
$\gamma$ in $G^{(1)}$
so that $\{(s(\gamma^\prime), r(\gamma^\prime)) \: : \: \gamma^\prime \in U\}$ 
is the graph of a
diffeomorphism $h \: : \:  V \rightarrow W$ between neighborhoods 
(in $G^{(0)}$) of $s(\gamma)$ and $r(\gamma)$.
In this way, $G$ gives rise to
a {\em pseudogroup} of diffeomorphisms of $G^{(0)}$. Conversely,
if ${\mathcal P}$ is a pseudogroup of diffeomorphisms of $G^{(0)}$ then
there is a corresponding \'etale groupoid $G$, with $G^{(1)}$ consisting of 
the germs of elements of ${\mathcal P}$. We say that an \'etale groupoid $G$
is {\em effective} when the germ of an element 
$\gamma \in G^{(1)}$ is trivial
if and only if $\gamma \in G^{(0)}$. An example of a noneffective groupoid
comes from a nontrivial discrete group $\Gamma$, with
$G^{(1)} \: = \: \Gamma$ and $G^{(0)} \: = \: \{e\}$.
Hereafter the \'etale groupoids will be assumed to be effective.

An \'etale groupoid is {\em Riemannian} if there is a Riemannian
metric $g$ on $G^{(0)}$ so that the germs of elements of $G^{(1)}$ act
isometrically. There is an obvious notion of the length of a smooth path
in $G$. There is a pseudometric $d$ on the orbit space
${\mathcal O}$, given by saying that for $x, y \in G^{(0)}$,
$d(O_x, O_y)$ is the infimum of the lengths of the smooth paths joining
$x$ to $y$. The Riemannian groupoid $G$
is {\em complete} if $({\mathcal O}, d)$
is a complete pseudometric space. Hereafter the Riemannian groupoids will
be assumed to be complete.
The {\em diameter} of $G$ is the
diameter of the pseudometric space $({\mathcal O}, d)$.
If $(G, O_x)$ is a pointed Riemannian groupoid then we write
$B_R(O_x) \: = \: \bigcup \{O_y \: : \: d(O_y, O_x) \: < \: R\} \subset
G^{(0)}$.
Two Riemannian groupoids are (pointed) 
{\em isometrically equivalent} if there is a (pointed) isometric equivalence
between them, as defined in terms of localizations.

Let
$J_1$ be the groupoid of $1$-jets of local diffeomorphisms of
$G^{(0)}$. With the natural topology on $J_1^{(1)}$,
it is a smooth non\'etale groupoid with $J_1^{(0)} \: = \:
G^{(0)}$. There is a continuous morphism $j_1 \: : \: G \rightarrow J_1$
which is injective,  as the germ of an isometry is determined by
its $1$-jet. Taking the closure of $j_1 \left( G^{(1)} 
\right)$ in $J_1^{(1)}$ gives a space that can be written as the embedding
$\overline{j}_1 \: : \: \overline{G} \rightarrow J_1$ of a unique
Riemannian groupoid
$\overline{G}$ with $\overline{G}^{(0)} \: = \: {G}^{(0)}$. 
It is called the {\em closure} of $G$. In effect, 
$\overline{G}^{(1)}$ is obtained by taking the closure of
$j_1 \left( {G}^{(1)} \right)$ and changing the topology to give it
the structure of an \'etale groupoid.

The orbits of $\overline{G}$ are closed submanifolds of
$G^{(0)}$. The orbit space of $\overline{G}$ is Hausdorff.
There is a locally constant sheaf ${\underline{\frak g}}$ on
$G^{(0)}$ of finite-dimensional Lie algebras, called the {\em structure
sheaf}, with the following properties : \\
1. ${\underline{\frak g}}$ is a $\overline{G}$-sheaf. \\
2. ${\underline{\frak g}}$ is a sheaf of germs of Killing vector fields. \\
3. The elements of $\overline{G}$ that are close to units in $J_1$ are 
germs of local isometries $\exp (\xi)$, where $\xi$ is a 
local section $\xi$ of ${\underline{\frak g}}$ that is close to zero.

Hereafter we assume that
the Riemannian groupoids are closed, unless otherwise stated.

\begin{example} \label{example5}
If $G \: = \: S^1 \rtimes \Z$ is the groupoid of Example \ref{example2},
with the generator of $\Z$ acting by an irrational rotation, then
$\overline{G} \: = \: S^1 \rtimes \SO(2)$ and ${\underline{\frak g}}$ is the
constant $\R$-sheaf on $S^1$.
\end{example}

\begin{example} \label{example6}
Let $M$ be a complete Riemannian manifold with sectional curvatures
between $-K$ and $K$, for some $K > 0$. Given
$r \: < \: \frac{\pi}{\sqrt{K}}$, for any $m \in M$ the exponential map
$\exp_m \: : \: T_mM \rightarrow M$ 
restricts to a local diffeomorphism from the $r$-ball 
$B_r^{(m)}(0) \subset T_mM$ to $B_r(m) \subset M$. Put the metric
$(\exp_m)^* g$ on $B_r^{(m)}(0)$.

Let $\{m_i\}_{i \in I}$ be points in $M$ so that $\{B_r(m_i)\}_{i \in I}$
covers $M$. Define a Riemannian groupoid $G$ with
$G^{(1)} \: = \: \bigsqcup_{i,j \in I} \{ (v_i, v_j) \in
B_r^{(m_i)}(0) \times B_r^{(m_j)}(0) \: : \:
\exp_{m_i}(v_i) \: = \: \exp_{m_j}(v_j) \}$ and  
$G^{(0)} \: = \: \bigsqcup_{i \in I} B_r^{(m_i)}(0)$ by
$r(v_i, v_j) \: = \: v_i$, $s(v_i, v_j) \: = \: v_j$ and
$(v_i, v_j) \cdot (v_j, v_k) \: = \: (v_i, v_k)$. 
Then $G$ is isometrically equivalent to the groupoid $M$.
\end{example}

\subsection{Convergence of Riemannian Groupoids}

\begin{definition} \label{defconv}
Let $\{(G_i, O_{x_i})\}_{i=1}^\infty$ be a sequence of pointed
$n$-dimensional 
Riemannian groupoids. Let $(G_\infty, O_{x_\infty})$ be a pointed
Riemannian groupoid. Let $J_1$ be the groupoid of 
$1$-jets of local diffeomorphisms of $G_\infty^{(0)}$.
We say that $\lim_{i \rightarrow \infty} (G_i, O_{x_i}) \: = \: 
(G_\infty, O_{x_\infty})$ in the pointed smooth topology if for all 
$R > 0$, \\
1.  There are pointed diffeomorphisms 
$\phi_{i, R} \: : \: B_{R}(O_{x_\infty}) \rightarrow B_{R}(O_{x_i})$, defined for
large $i$, from the 
pointed  $R$-ball in $G_\infty^{(0)}$ to the pointed $R$-ball in
$G_i^{(0)}$, so that $\lim_{i \rightarrow \infty} 
\phi_{i, R}^* \: g_i \big|_{B_{R}(O_{x_i})} \: = \:
g_\infty \big|_{B_{R}(O_{x_\infty})}$. \\
2. After conjugating by $\phi_{i, R}$, the images of
$s_i^{-1}(\overline{B_{R/2}(O_{x_i})}) \cap 
r_i^{-1}(\overline{B_{R/2}(O_{x_i})})$ in 
$J_1$ converge in the Hausdorff sense to the image of 
$s_\infty^{-1}(\overline{B_{R/2}(O_{x_\infty})}) \cap 
r_\infty^{-1}(\overline{B_{R/2}(O_{x_\infty})})$ in $J_1$.
\end{definition}

This definition is similar to 
\cite[Definition A.4]{Petrunin-Tuschmann (1999)}. The paper
\cite{Petrunin-Tuschmann (1999)} considers Lipschitz convergence instead
of smooth convergence.

We will allow ourselves to freely replace a Riemannian groupoid by an
isometrically equivalent one,
without saying so explicitly.

\begin{proposition} \label{limits}
Let $\{(M_i, p_i)\}_{i=1}^\infty$ be a sequence of pointed complete
$n$-dimensional Riemannian
manifolds. Suppose that for each $a \in \Z^{\ge 0}$ and $R > 0$, there is
some $K_{a, R} < \infty$ so that for all $i$, one has
$\parallel \nabla^a\Riem(M_i) \parallel_\infty \: \le \: K_{a,R}$ on
$B_R(p_i)$. Then there is a subsequence of $\{(M_i, p_i)\}_{i=1}^\infty$
that converges to some pointed $n$-dimensional
Riemannian groupoid $(G_\infty, O_{x_\infty})$
in the pointed smooth topology. 
\end{proposition}
\begin{proof}
Put $r(j) \: = \: \frac{\pi}{2\sqrt{K_{0,2j}}}$.
With reference to Example \ref{example6}, for each $j \in \Z^+$ there
is a number $N_j$ so that
we can find points
$\{x_{i,j,k}\}_{k=1}^{N_j}$ in $B_j(p_i) - B_{j-1}(p_i)$ with the
property that
$\bigcup_{j=1}^\infty \bigcup_{k=1}^{N_j} B_{r(j)}(x_{i,j,k})$
covers $M_i$. As in Example \ref{example6}, we form the corresponding
Riemannian groupoid $G_i$ with $G_i^{(0)} \: = \: 
\bigsqcup_{j=1}^\infty \bigsqcup_{k=1}^{N_j} B^{(x_{i,j,k})}_{r(j)}(0)$. It is
isometrically equivalent to the Riemannian groupoid $M_i$.
After passing to a subsequence,
we may assume that $\{G_i^{(0)}\}_{i=1}^\infty$ converges smoothly
to some $G_\infty^{(0)}$ in the sense of Definition \ref{defconv}.1.
The rest of the argument is basically the same as in
\cite[Pf. of Theorem 0.5]{Fukaya (1988)}, which in turn uses ideas from
\cite[Chapitres 8C and 8D]{Gromov (1981b)}.
Namely, after passing to a further subsequence, we can construct
$G_\infty^{(1)}$ as a pointed Hausdorff limit, in the sense of
Definition \ref{defconv}.2, of the images of
$G_i^{(1)}$ in $J_1^{(1)}$, where $J_1$ is the groupoid of $1$-jets
of local diffeomorphisms of $G_\infty^{(0)}$. (Because of
the convergence of the metrics in the sense of Definition \ref{defconv}.1,
the images of $G_i^{(1)}$ in $J_1^{(1)}$ come closer and closer
to taking value in the $1$-jets of local isometries.)

\end{proof}

From the construction of $G_\infty$, there is a
$G_\infty$-invariant sheaf ${\underline{\frak g}}$ on
$G_\infty^{(0)}$ consisting of local Killing vector fields that 
generate the collapsing directions.  From 
\cite[Section 4]{Fukaya (1988)}, 
${\underline{\frak g}}$ is a sheaf of nilpotent Lie algebras.

The orbit space of $G_\infty$ is the same as the pointed Gromov-Hausdorff
limit of the convergent subsequence of $\{(M_i, p_i)\}_{i=1}^\infty$.

Proposition \ref{limits} is essentially equivalent to
\cite[Theorem A.5(i)]{Petrunin-Tuschmann (1999)}, except that we
use smooth convergence instead of Lipschitz convergence. 
Ricci flow will provide the needed smoothness.

\begin{example}
Let $M_i$ be the circle of radius $i^{-1}$. Then
$\lim_{i \rightarrow \infty} \: M_i$ is the cross-product 
groupoid $\R \rtimes \R$ of Example \ref{example2}.
\end{example}

\begin{example}
As in Example \ref{sl2example},
let $\Sigma$ be a closed hyperbolic surface and let
$S^1 \Sigma$ be its unit sphere bundle.  Let
$M_i$ be $S^1 \Sigma$ equipped with $\frac{1}{i}$ times the time-$i$ Ricci
flow metric.  Then $\lim_{i \rightarrow \infty} (M_i, p)$ is
the cross-product groupoid $(\R \times \Sigma) \rtimes \R$. 
\end{example}

\subsection{Convergence of Ricci Flows on \'Etale Groupoids}

Let $G$ be a complete Riemannian groupoid.
We can construct its Ricci tensor $\Ric(g)$ as a symmetric 
covariant $2$-tensor field
on $G^{(0)}$ which is invariant in the sense that it is preserved by
the germs of elements of $G^{(1)}$.

Let $\{g(t)\}$ be a smooth $1$-parameter family of Riemannian metrics on
the \'etale groupoid
$G$, where smoothness in $t$ can be checked locally on $G^{(0)}$.
Then $\{g(t)\}$ satisfies the Ricci flow equation if
$\frac{\partial g}{\partial t} \: = \: - \: 2 \Ric$. We say that it is
an {\em expanding soliton} if the flow is defined for $t \in (0, \infty)$ and
$g(t)$ is isometrically equivalent to
$t \: g(1)$ for all $t \in (0, \infty)$.

We now state an analog of
\cite[Theorem 1.2]{Hamilton (1995)}, except without the assumption
of a lower bound on the injectivity radius.
We define convergence of a sequence of Ricci flow solutions as
in Section \ref{def2},  except that now we allow the limit to be a
Ricci flow on an \'etale groupoid.

\begin{theorem} \label{RFlimits}
Let $\{(M_i, p_i, g_i(\cdot))\}_{i=1}^\infty$ be a sequence of
Ricci flow solutions on pointed $n$-dimensional
manifolds $(M_i, p_i)$. We assume that 
there are numbers $-\infty \: \le A \: < \: 0$ and $0 \: < \:\Omega
\: \le \: \infty$ so that \\
1. The Ricci flow solution $(M_i, p_i, g_i(\cdot))$ is defined on the
time interval $(A, \Omega)$. \\
2. For each $t \in (A,\Omega)$, $g_i(t)$ is a complete Riemannian metric
on $M_i$. \\
3. For each closed interval $I  \subset (A, \Omega)$
there is some $K_{I} \: < \: \infty$ so that $|\Riem(g_i)(x, t)| \: \le \:
K_{I}$ for all $x \in M_i$ and $t \in I$. \\

Then after passing to a subsequence,
the Ricci flow solutions $g_i (\cdot)$ converge
smoothly to 
a Ricci flow solution $g_\infty(\cdot)$ on a pointed $n$-dimensional
\'etale groupoid 
$\left( G_\infty, O_{x_\infty} \right)$, defined again for 
$t \in (A, \Omega)$.
\end{theorem}

\begin{proof}
The proof is virtually the same as that of
\cite[Theorem 1.2]{Hamilton (1995)}. From local derivative estimates,
the pointed Riemannian manifolds 
$\{(M_i, p_i, g_i(0))\}_{i=1}^\infty$ satisfy the assumptions of
Proposition \ref{limits}. Then after passing to a subsequence, we can assume
that $\{(M_i, p_i, g_i(0))\}_{i=1}^\infty$ converges to a pointed
Riemannian groupoid $\left( G_\infty, O_{x_\infty}, g_\infty(0) \right)$.
In terms of the proof of Proposition \ref{limits}, we have pointed
time-$0$ convergence of the Riemannian groupoids $G_i$ to $G_\infty$. 
The remaining step is to
get convergence on the whole time interval $(A, \Omega)$, after
passing to a further subsequence. The argument for this is as in
\cite[Section 2]{Hamilton (1995)}. 
Namely, for any $R > 0$, if $i$ is sufficiently
large then one can use $\phi_{i,R}$ to 
transfer the time-$t$ metric $g_i(t)$ on
$B_R \left( O_{x_i} \right) \subset G_i^{(0)}$ to the time-$t$ metric 
$\phi_{i,R}^* g_i(t)$ on the time-zero set
$B_R \left( O_{x_\infty} \right) \subset G_\infty^{(0)}$.
As in \cite[Lemma 2.4]{Hamilton (1995)}, for any closed subinterval $I
\subset (A, \Omega)$ one has uniform bounds on the norms
(with respect to $g_\infty$) of the
spatial and temporal derivatives of 
$\phi_{i,R}^* g_i(\cdot)$ on $I \times B_R \left( O_{x_\infty} \right)$.
Then after passing to a further subsequence, one obtains a limiting
metric $g_\infty(\cdot)$ on $G_\infty$ which will necessarily satisfy
the Ricci flow equation.
\end{proof}

\begin{remark}
If one considers individual tangent spaces $T_{p_i} M_i$ instead of groupoids
then an analog of
\cite[Theorem 1.2]{Hamilton (1995)} was proven in \cite{Glickenstein (2003)}.
The result of \cite{Glickenstein (2003)} was used in 
\cite{Chow-Glickenstein-Lu (2003)} to study three-dimensional
type-II Ricci flow solutions.
\end{remark}

\begin{corollary}
Let $(M, p, g(\cdot))$ be a type-III Ricci flow solution. 
Then for any
sequence $s_i \rightarrow \infty$, there are a subsequence (which we
relabel as $\{s_i \}_{i=1}^\infty$) and a Ricci flow solution 
$\left( G_\infty, O_{x_\infty}, g_\infty(\cdot) \right)$, defined for
$t \in (0, \infty)$, so that
$\lim_{i \rightarrow \infty} (M, p, g_{s_i}(\cdot)) \: = \: 
\left( G_\infty, O_{x_\infty}, g_\infty(\cdot) \right)$.
\end{corollary}

\begin{corollary} \label{compactness}

Given $K > 0$, the space of pointed $n$-dimensional Ricci flow solutions with
$\sup_{t \in (0, \infty)} \: t \: \parallel \Riem(g_t) \parallel_\infty
\: \le \: K$ is relatively compact among Ricci flows on
pointed $n$-dimensional \'etale
groupoids, defined for $t \in (0, \infty)$.
\end{corollary}

\subsection{Compactification of type-III Ricci flow solutions} \label{discussion}

With reference to Corollary \ref{compactness}, 
let ${\mathcal S}_{n,K}$ denote the closure of the pointed
$n$-dimensional  Ricci flow
solutions on manifolds
with $\sup_{t \in (0, \infty)} \: t \: \parallel \Riem(g_t) \parallel_\infty
\: \le \: K$. Given $g(\cdot) \in {\mathcal S}_{n,K}$,
there is a rescaled Ricci flow solution $g_s(\cdot) \in 
{\mathcal S}_{n,K}$ given by
$g_s(t) \: = \: s^{-1} \: g(st)$. This means that there is
an $\R^+$-action on ${\mathcal S}_{n,K}$ where $s \in \R^+$
sends $g$ to $g_s$. 
Understanding the long-time behavior of a type-III 
Ricci flow solution $g$ boils down to understanding the dynamics of
its orbit in the compact set ${\mathcal S}_{n,K}$. 
We make some elementary comments about ${\mathcal S}_{n,K}$. 

First, a Ricci flow in the boundary of ${\mathcal S}_{n,K}$ necessarily has a
collapsing structure, i.e. a nontrivial sheaf of nilpotent Lie algebras that
act as local Killing vector fields on the space of units.
For simple examples of Ricci flows in the boundary of ${\mathcal S}_{n,K}$,
consider first the Ricci flow of a generic flat metric on the $j$-torus.
Its rescaling limit is the constant Ricci flow on the cross-product groupoid
$\R^j \rtimes \R^j$, where the first $\R^j$ factor has a flat metric $g_0$.
Now if $(\widehat{M}, \widehat{g}(\cdot))$ is any pointed Ricci flow on an 
$(n-j)$-dimensional
manifold $\widehat{M}$ with 
$\sup_{t \in (0, \infty)} \: t \: \parallel \Riem(\widehat{g}_t) \parallel_\infty
\: \le \: K$ then the product flow $\widehat{g}(\cdot) + g_0$
on the groupoid $(\widehat{M} \times \R^j) \rtimes \R^j$ is an element
of the boundary of ${\mathcal S}_{n,K}$, as it is a limit of Ricci flows
on $\widehat{M} \times T^j$. The $\R^+$-action commutes 
with the inclusion ${\mathcal S}_{n-j,K} \rightarrow 
\partial {\mathcal S}_{n,K}$.

Now let ${\mathcal E}_{n,K}$ be the $n$-dimensional pointed 
Einstein metrics $g$ with
Einstein constant $- \: \frac12$ and $\parallel \Riem(g) \parallel_\infty
\: \le \: K$,  modulo pointed diffeomorphisms. 
We can identify ${\mathcal E}_{n,K}$ with a set of
Ricci flow solutions, with the Einstein metric $g$ being the
time-$1$ metric of the solution. In this way there is an inclusion
${\mathcal E}_{n,K} \subset {\mathcal S}_{n,K}$. We remark
that by moving the basepoint on a finite-volume noncompact
manifold with constant sectional curvature $- \: \frac{1}{2(n-1)}$,
we obtain examples where $\overline{{\mathcal E}_{n,K}}$ intersects
$\partial {\mathcal S}_{n,K}$.

Consider a pointed $n$-dimensional compact Ricci flow solution 
$(M, p, g(\cdot))$ with $\sup_{t \in (0, \infty)} \: t \: \parallel \Riem(g_t) \parallel_\infty
\: \le \: K$. Suppose that for large $s$, the corresponding orbit 
$\{g_s(\cdot)\}$ in ${\mathcal S}_{n,K}$ stays away from the
boundary. This is equivalent to saying that 
$\liminf_{t \rightarrow \infty} t^{-\frac12} \: \inj_{g(t)}(p) \: > \: 0$.
If there is a pointed limit $\lim_{j \rightarrow \infty} g_{s_j}(\cdot)$, with 
$\{s_j\}_{j=1}^\infty$ a sequence tending to infinity, then the limit is
a Ricci flow on a finite-volume Einstein manifold with negative Einstein constant
\cite[Section 7]{Hamilton (1999)},\cite[Section 7.1]{Perelman (2003a)}.
(The proof in \cite[Section 7]{Hamilton (1999)} is for the
normalized Ricci flow while the proof in 
\cite[Section 7.1]{Perelman (2003a)}
is for the unnormalized Ricci flow.  Both proofs are for $n = 3$ but
extend to general $n$.) It follows that as $s \rightarrow
\infty$, the orbit $\{g_s(\cdot)\}$ approaches ${\mathcal E}_{n,K}$. However,
it does not immediately follow that the orbit approaches a fixed-point.

To understand the asymptotics of the orbits that do not stay away from
the boundary of ${\mathcal S}_{n,K}$, 
a first question is whether the orbit approaches the
boundary, i.e. whether Ricci flow favors the formation of continuous symmetries.

Independent of this question, one can ask about
the dynamics of the $\R^+$-action on the boundary.
As the boundary elements have a collapsing structure, in principle one can use
this to help analyze the Ricci flow equations. An example is 
Proposition \ref{helimits}, where a monotonicity formula was used.
An overall question is whether
the orbits of the $\R^+$-action on ${\mathcal S}_{n,K}$ approach fixed
points, i.e. expanding soliton solutions.

In the case of a finite-volume locally homogeneous three-manifold,
the results of Section \ref{homogeneous} imply that the
orbit approaches a fixed-point $g_\infty(\cdot)$.
The next proposition makes this explicit. In the statement of the
proposition
we allow the fundamental group $\Gamma$ to be a discrete subgroup of
a {\em maximal} element $G$ among groups of diffeomorphisms of the universal
cover that act transitively with compact isotropy group, and we take the
Riemannian metric to be $G$-invariant;
see \cite[\S 4 and \S 5]{Scott (1983)} and
\cite[Chapter 3]{Thurston (1997)} for the description
of such groups. 

\begin{proposition}
Let $(M, p, g(\cdot))$ be a finite-volume pointed
locally homogeneous three-dimensional
Ricci flow solution that exists for all $t \in (0, \infty)$.  
Then $\lim_{s \rightarrow \infty} (M, p, g_s(\cdot))$
exists and is an expanding soliton on a pointed three-dimensional
\'etale groupoid $G_\infty$.
Put $\Gamma \: = \: \pi_1(M, p)$. 
The groupoid $G_\infty$ and its metric $g_\infty(t)$ are given as follows. \\
1. If $(M, g(0))$ has constant negative curvature then
$G_\infty$ is the cross-product groupoid $H^3 \rtimes \Gamma$
(which is equivalent to $M$)
and $g_\infty(t)$ has constant sectional curvature $- \: \frac{1}{4t}$.\\
2. If $(M, g(0))$ has $\R^3$-geometry, there is a homomorphism
$\alpha \: : \: \Isom(\R^3) \rightarrow \Isom(\R^3)/\R^3$. 
(Here $\R^3$ is the translation subgroup of $\Isom(\R^3)$. The quotient
$\Isom(\R^3)/\R^3$ is isomorphic to $O(3)$.) Put
$\Gamma_\R \: = \: \alpha^{-1}(\alpha(\Gamma))$.
Then $G_\infty$ is the cross-product groupoid
$\R^3 \rtimes \Gamma_\R$
and $g_\infty(t)$ is the constant flat metric. \\
3. If $(M, g(0))$ has $\Sol$-geometry, there is a homomorphism
$\alpha \: : \: \Isom(\Sol) \rightarrow \Isom(\Sol)/\R^2$. 
(Here $\R^2$ is a normal subgroup of $\Sol$, which is a normal
subgroup of $\Isom(\Sol)$. The quotient $\Isom(\Sol)/\R^2$ has
$\R$ as a normal subgroup of index $8$.) Put
$\Gamma_\R \: = \: \alpha^{-1}(\alpha(\Gamma))$. Then
$G_\infty$ is the cross-product groupoid $\Sol \rtimes \Gamma_\R$, with
$g_\infty(t)$ given by (\ref{solmetric}).\\
4. If $(M, g(0))$ has $\Nil$-geometry, there is a homomorphism
$\alpha \: : \: \Isom(\Nil) \rightarrow \Isom(\Nil)/\Nil$. 
(Here $\Nil \subset \Isom(\Nil)$ acts by left multiplication.
The quotient $ \Isom(\Nil)/\Nil$ is isomorphic to $O(2)$.) Put
$\Gamma_\R \: = \: \alpha^{-1}(\alpha(\Gamma))$. Then
$G_\infty$ is the cross-product groupoid $\Nil \rtimes \Gamma_\R$, with
$g_\infty(t)$ given by (\ref{nilmetric}).\\
5. If $(M, g(0))$ has $(\R \times H^2)$-geometry, there is a 
homomorphism $\alpha \: : \: \Isom \left( 
\R \times H^2 \right) \rightarrow \Isom \left( 
\R \times H^2 \right)/\R$. 
(Here $\Isom \left( \R \times H^2 \right)/\R$ is isomorphic to 
$\Z_2 \times \Isom(H^2)$.)
Put
$\Gamma_\R \: = \: \alpha^{-1}(\alpha(\Gamma))$.
Then
$G_\infty$ is the cross-product groupoid
$(\R \times H^2) \rtimes \Gamma_\R$ 
and $g_\infty(t) \: = \: g_\R \: + \: g_{H^2}(t)$,
where $g_{H^2}(t)$ has constant sectional curvature $- \: \frac{1}{2t}$. \\
6. If $(M, g(0))$ has $\widetilde{\SL_2(\R)}$-geometry, there is a 
homomorphism $\alpha \: : \: \Isom \left( 
\widetilde{\SL_2(\R)} \right) \rightarrow \Isom \left( 
\widetilde{\SL_2(\R)} \right)/\R$.
(Here $\Isom \left( 
\widetilde{\SL_2(\R)} \right)/\R$ is isomorphic to $\Isom(H^2)$.)
Then $G_\infty$ is the cross-product groupoid
$(\R \times H^2) \rtimes 
(\R \widetilde{\times} \alpha({\Gamma}))$ 
and $g_\infty(t) \: = \: g_\R \: + \: g_{H^2}(t)$,
where $g_{H^2}(t)$ has constant sectional curvature $- \: \frac{1}{2t}$.
In writing $\R \widetilde{\times} \alpha({\Gamma})$, the group
$\alpha(\Gamma) \subset \Isom(H^2)$ acts linearly on $\R$ via the 
orientation homomorphism $\alpha(\Gamma) \rightarrow \Z_2$.
\end{proposition}
\begin{proof}
This follows from the results of Section \ref{3d}. We just give two examples.

In case 4, suppose for simplicity that $\Gamma \subset \Nil$ and
$M \: = \: \Gamma \backslash \Nil$ is compact.
Fix a basepoint $p \in M$. Following the proof of Theorem \ref{RFlimits},
we will first construct a limiting time-$1$ Riemannian groupoid 
$(G_\infty, O_{x_\infty}) \: = \: \lim_{s \rightarrow \infty} \left(M, p, g_s(1)
\right)$.  From (\ref{nilflow}), we have
$\lim_{s \rightarrow \infty} \: \diam(M, g_s(1)) \: = \: 0$.
For any $r > 0$, if $s$ is sufficiently large then the exponential map
$\exp_{p}(s) \: : \: B_r^{(p)}(0) \rightarrow M$ for $M$, with the metric $g_s(1)$, 
is a surjective local
diffeomorphism. From the calculations in Section \ref{nil3}, after performing
appropriate diffeomorphisms the pullback metrics $\exp_p(s)^* g_s(1)$ will
approach $\exp_{p_\infty}^* g_\infty(1)$, where $g_\infty(1)$ is the metric of
(\ref{nilmetric}) and 
$\exp_{p_\infty} \; : \: B_r^{(p_\infty)}(0) \rightarrow M_\infty$ is the corresponding
exponential map on the $r$-ball in the tangent space
at the basepoint $p_\infty \in M_\infty$. 
The limit groupoid $G_\infty$ consists of germs of
local isometries of $B_r^{(p_\infty)}(0)$, the latter being equipped with the
Riemannian metric $\exp_{p_\infty}^* g_\infty(1)$. This groupoid
is isometrically equivalent
to the cross-product groupoid $\Nil \rtimes \Nil$, with metric $g_\infty(1)$.
The extension to times $t$ other than $1$ gives the Ricci flow on
the cross-product groupoid $\Nil \rtimes \Nil$ with metric $g_\infty(\cdot)$.

In case 6, suppose that $M$ is the pointed unit sphere bundle $S^1 \Sigma$ of an
oriented hyperbolic
surface $\Sigma$. The pointed Gromov-Hausdorff limit $\lim_{s \rightarrow \infty}
(M, p, g_s(1))$ is the surface $\Sigma$ with the metric rescaled to have 
sectional curvature $- \: \frac12$. We choose points $x_{j,k}$ as in the proof of
Proposition \ref{limits}; we can choose them to be independent of the
parameter $s$ (which replaces $i$). In the limit we obtain a groupoid
$(G_\infty, O_{x_\infty})$ whose unit space 
$G_\infty^{(0)}$consists of a disjoint union of
manifolds, each isometrically equivalent to a domain in
$\R \times H^2$; see Example \ref{sl2example}. After making the remaining
identifications, this groupoid is
isometrically equivalent to
$(\R \times \Sigma) \rtimes \R$, equipped with the product metric
$g_\R \; + \: g_{\Sigma}(1)$. The extension to times $t$ other than $1$ gives the
Ricci flow on $(\R \times \Sigma) \rtimes \R$ with metric $g_\R \; + \: g_{\Sigma}(t)$,
where $g_{\Sigma}(t)$ has constant sectional curvature $- \: \frac{1}{2t}$.
This, in turn, is isometrically equivalent to the cross-product groupoid
$(\R \times H^2) \rtimes (\R \times \alpha(\Gamma))$ 
with the metric $g_\R \; + \: g_{H^2}(t)$.
\end{proof}

\bibliographystyle{acm}

\end{document}